\documentclass[psamsfonts]{amsart}
\usepackage[utf8]{inputenc}
\usepackage{amsfonts}
\usepackage{amsmath}
\usepackage{xcolor}
\usepackage{amsthm}
\usepackage{pdflscape}
\usepackage{pgfplots}
\usepackage{mathrsfs}

\newtheorem{thm}{Theorem}[section]
\newtheorem*{thm*}{Theorem}
\newtheorem{cor}[thm]{Corollary}
\newtheorem{prop}[thm]{Proposition}
\newtheorem{lem}[thm]{Lemma}

\newtheorem{quest}[thm]{Question}

\newtheorem{notation}[thm]{Notation}

\theoremstyle{definition}
\newtheorem{defn}[thm]{Definition}

\newtheorem{ex}[thm]{Example}

\theoremstyle{remark}
\newtheorem{rmk}[thm]{Remark}

\usepackage{amsmath}	
\usepackage{amssymb}		
\usepackage{amsthm}		
\usepackage{setspace}	
\usepackage{array} 
\usepackage{float}
\usepackage[mathscr]{euscript} 
 \let\mathscr\relax
\usepackage[scr]{rsfso}

\newcommand{\m}{\mathfrak{m}}

\newcommand{\la}{\langle}
\newcommand{\ra}{\rangle}

\newcommand{\Hom}{\text{Hom}}

\DeclareMathOperator{\coker}{coker}
\newcolumntype{P}[1]{>{\centering\arraybackslash}p{#1}}		            
\usepackage[top=1in,bottom=1in,left=1.25in,right=1in]{geometry}	

\newtheorem*{lemma*}{lemma}

\newtheorem{ques}[thm]{Question}
\usepackage{graphicx}

\usepackage[all]{xy}
\newcommand{\im}{\hspace{1mm}\text{im}\hspace{1mm}}

\newcommand{\onto}{\twoheadrightarrow}

\DeclareMathOperator{\Der}{Der}
\DeclareMathOperator{\grade}{grade}
\DeclareMathOperator{\height}{ht\hspace{-0.5mm}}

\usepackage{tikz-cd}
\usetikzlibrary{cd}
\usepackage{adjustbox}
\usepackage{caption}  
\usepackage{subcaption}  

\usepackage{xurl}
\usepackage[colorlinks]{hyperref}
\hypersetup{
 colorlinks=true,
 linkcolor=blue,
 filecolor=magenta,
 urlcolor=cyan,
 breaklinks=true,
}

\usepackage[maxbibnames=9]{biblatex} 
\makeatletter
\let\c@equation\c@thm
\makeatother
\numberwithin{equation}{section}

\title{Resolving the Module of Derivations on an $n \times (n+1)$ Determinantal Ring}

\author{Henry Potts-Rubin}
\address{Mathematics Department, Syracuse University, Syracuse, NY 13244 U.S.A.}
\email{\href{mailto:hpottsru@syr.edu}{hpottsru@syr.edu}}
\keywords{resolution, derivations, determinantal ring}

\begin{document}

\maketitle
\vspace{-10mm}
\begin{abstract}
    We use the construction of the relative bar resolution via differential graded structures to obtain the minimal graded free resolution of $\Der_{R \mid k}$, where $R$ is a determinantal ring defined by the maximal minors of an $n \times (n+1)$ generic matrix and $k$ is its coefficient field.  Along the way, we compute an explicit action of the Hilbert-Burch differential graded algebra on a differential graded module resolving the cokernel of the Jacobian matrix whose kernel is $\Der_{R \mid k}$.  As a consequence of the minimality of the resulting relative bar resolution, we get a minimal generating set for $\Der_{R \mid k}$ as an $R$-module, which, while already known, has not been obtained via our methods.  
\end{abstract}


\section{Introduction}

This paper is concerned with showing that a truncation of a relative bar resolution minimally resolves the module of $k$-linear derivations $\Der_{R \mid k}$ on a determinantal ring $R$ defined by the maximal minors of an $n \times (n+1)$ generic matrix $X$.   The minimality of the resolution yields a minimal generating set for $\Der_{R \mid k}$ over $R$, a result which may be found in \cite{bv} or \cite{CB}, each produced by a starkly different methods.     


While minimal generators of modules of derivations on determinantal rings are known, the minimal graded free resolution of $\Der_{R \mid k}$ over $R$ has not been computed.  The first order differential operators $D^1_{R \mid k}$ on $R$ decompose as $D^1_{R \mid k}\cong\Der_{R \mid k} \oplus R$, and so after obtaining the resolution of $\Der_{R \mid k}$, the resolution of $D^1_{R \mid k}$ is immediate.  Known resolutions of differential operators are few $-$ see \cite{rc} for the case of low-order differential operators on an isolated hypersurface singularity and \cite{yu} for the case of derivations on a generic hyperplane arrangement, as two examples.  Differential operators need not behave the same in positive characteristic as in characteristic zero, so it seems worthwhile to point out that our approach to the minimal graded $R$-free resolution of $\Der_{R \mid k}$ is characteristic-free.  

The $R$-module $\Der_{R \mid k}$ is isomorphic to the kernel of $J^T_R$, where $J^T_R$ is the transpose of the Jacobian matrix of the maximal minors of $X$ tensored down to $R$.  That is, the following sequence of $R$-modules is exact:
\[
0 \to \Der_{R \mid k} \to R^{n(n+1)} \xrightarrow{J^T_R} R^{n+1} \to \coker J^T_R \to 0. 
\]
Consequently, to obtain the minimal graded $R$-free resolution of $\Der_{R \mid k}$, we may construct the minimal graded $R$-free resolution of $\coker J^T_R$ and truncate it.  

Our main result is the following.  For the proof, we build the appropriate differential graded (dg) structures to use the construction of the relative bar resolution developed by Iyengar in \cite{Iyengar}.      

\begin{thm*}[Theorem \ref{mainresult}]
    Let $\mathcal{A}$ be the minimal graded $Q$-free resolution of $R$ and $\mathcal{U}$ be the minimal graded $Q$-free resolution of $\coker J^T_R$.  The minimal graded $R$-free resolution of $\Der_{R \mid k}$ is the truncated relative bar resolution
    \[
    \dots \longrightarrow R^{n(n+1)(n+2)}\xrightarrow{\small\begin{bmatrix}
        J^T_R & \varphi_2\\
        M_{1,1} & M_{2,0}
    \end{bmatrix}}R^{2n(n+1)}\xrightarrow{\small\begin{bmatrix}
        \partial_2 & M_{1,0}
    \end{bmatrix}} \Der_{R \mid k} \longrightarrow 0,
    \]
    where $M_{i,j}$ comes from the action $\mathcal{A}_i\cdot \mathcal{U}_j$ of the dg algebra $\mathcal{A}$ on the dg $\mathcal{A}$-module $\mathcal{U}$ and $\varphi_2$ and $\partial_2$ come from the differentials of $\mathcal{A}$ and $\mathcal{U}$, respectively.  
\end{thm*}

It is not always the case that this construction produces a minimal resolution, but when the ring and module involved are Golod and the resolutions of them used in the construction are minimal, it does.  Our minimal free resolutions are graded, and we refrain from saying so from here on.  We do not indicate shifts.  In our setting, $\mathcal{A}$ is the minimal $Q$-free resolution of $R$, which is Golod, and $\mathcal{U}$ is the minimal $Q$-free resolution of $\coker J^T_R$, which is also Golod.  The resolution $\mathcal{A}$ is known to admit a dg algebra structure, and we show that $\mathcal{U}$ admits a dg $\mathcal{A}$-module structure.  The relative bar resolution then yields the minimal $R$-free resolution of $\coker J^T_R$, in which the first differential is $J^T_R$.  Truncating the relative bar resolution results in the minimal $R$-free resolution of $\Der_{R \mid k}$.  

The outline of sections is as follows.  Section \ref{background} introduces the necessary background, including the dg algebra structure on the Hilbert-Burch complex $\mathcal{A}$, and the specifics of this structure prove that $R$ is Golod.  Section \ref{lemmas}, and in particular Lemma \ref{syzarg}, discusses the relationship between minors and partial derivatives, to be used in several subsequent proofs.  In Section \ref{rescok}, we give the minimal $Q$-free resolution $\mathcal{U}$ of $\coker J^T_R$, where $J^T_R$ is the transpose of the Jacobian of the maximal minors of $X$ tensored down to $R$.  In Section \ref{dgastructure}, we compute an explicit dg algebra action of $\mathcal{A}$ on $\mathcal{U}$ which shows that $\coker J^T_R$ is Golod, so that in Section \ref{resder} we may use the relative bar resolution via dg structures to obtain the minimal $R$-free resolution of $\Der_{R \mid k}$. 


\section{Background}\label{background}

In this section, we recall the background necessary to talk about the ingredients going into the main result, Theorem \ref{mainresult}.  Throughout, let $k$ be a field and $Q=k[x_{i,j} \mid 1 \leq i \leq n,\,1 \leq j \leq n+1]$, where $n$ is at least 2.  Take $X$ to be the $n \times (n+1)$ matrix of variables with $(i,j)$-entry $x_{i,j}$.  For $r \in \{1,\ldots,n+1\}$, denote by $F_r$ the maximal minor of $X$ obtained by deleting the $r$th column.  For lists of positive integers $\alpha$ and $\beta$, write $X_\alpha^\beta$ for the matrix obtained by deleting the columns with labels in $\alpha$ and rows with labels in $\beta$, e.g., $F_r=\det X_r$.   

Write $I_n$ for the ideal of $Q$ generated by the $n \times n$ (i.e., maximal) minors of $X$, i.e., $I_n=(F_1,F_2,\ldots,F_{n+1})$, and set $R=Q/I_n$.  

\begin{defn}
The ring $R$ is called a \textit{determinantal ring}, and $I_n$ is called a \textit{determinantal ideal}.  
\end{defn}

More generally, determinantal rings are defined by the $k \times k$ minors of an $n \times m$ generic matrix; we simply focus on the case when $k=n$ and $m=n+1$ (and note that our methods do not easily generalize).  Much is known about determinantal rings and ideals $-$ see, for example, \cite{bv}.  One fact used in this work is that the height of $I_n$ is 2 (see, for example, \cite{bh}, Theorem 7.3.1).  Another relates to the minimal free resolution of $R$ over $Q$, but before stating it, we introduce differential graded algebras and differential graded modules.

\begin{defn}\label{dgadef}
A \textit{differential graded (dg) algebra} over a ring $S$ is a complex $(A,\partial^A)$ of free $S$-modules equipped with a unitary, associative multiplication $A \otimes_S A \to A$ satisfying
\begin{itemize}
    \item[(i)] $A_iA_j \subseteq A_{i+j}$,
    \item[(ii)] $\partial^A(a_ia_j)=\partial^A(a_i)a_j+(-1)^{i}a_i\partial^A(a_j)$,
    \item[(iii)] $a_ia_j=(-1)^{ij}a_ja_i$,
    \item[(iv)] $a_i^2=0 \text{ if $i$ is odd}$,
\end{itemize}
where $a_\ell \in A_\ell$.   
\end{defn}

\begin{defn}\label{dgmoddef}
Let $(A,\partial^A)$ be a dg algebra over a ring $S$.  A \textit{differential graded (dg) $A$-module} is a complex $(U,\partial^U)$ of free $S$-modules together with a unitary, associative multiplication $A \otimes_S U \to U$ satisfying
\begin{itemize}
    \item[(i)] $A_iU_j \subseteq U_{i+j}$,
    \item[(ii)] $\partial^U(a_iu_j)=\partial^A(a_i)u_j+(-1)^{i}a_i\partial^U(u_j)$,
\end{itemize}
where $a_i \in A_i$ and $u_j \in U_j$.  
\end{defn}  

A natural question is whether the minimal free resolution of a given module has the structure of a dg algebra.  For ideals of height 2 in a Cohen-Macaulay ring, the answer is given by the Hilbert-Burch complex, shown to be a resolution in the polynomial ring case by Hilbert in 1890 \cite{hilbert} and more generally by Burch in 1968 \cite{burch}.  The dg algebra structure is due to Herzog \cite{herzog}.  It is well-known that the height of $I_n$ is 2 (see, for example, \cite{bh}, Theorem 7.3.1), and so we state the next proposition in terms of $R$ and $Q$.        

\begin{prop}\cite{burch,herzog,hilbert}\label{EN}
    The minimal $Q$-free resolution of $R$ is the Hilbert-Burch complex
    \[
    0 \longrightarrow Q^n \xrightarrow{\varphi_2} Q^{n+1} \xrightarrow{\varphi_1} Q \longrightarrow R \longrightarrow 0 \tag{$\mathcal{A}$},
    \]
    where $\varphi_1=\begin{bmatrix}
        F_1&-F_2&\ldots&(-1)^{i+1}F_i&\ldots& (-1)^{n+2}F_{n+1}
    \end{bmatrix}$ and $\varphi_2=X^T$.  Furthermore, the resolution $\mathcal{A}$ is a dg algebra under the multiplication
    \begin{align*}
    e_i \cdot e_j=-(e_j \cdot e_i)&=\begin{cases}
   \displaystyle\sum_{k=1}^n(-1)^{i+j+k}\det X^k_{i,j}T_k, & \text{if } i < j,\\\\
        0, & \text{if } i=j,
    \end{cases}\\\\
    e_i \cdot T_j&=0,\\\\
    T_i \cdot T_j&=0,
    \end{align*}
    where $\{1\}$ is a basis for $\mathcal{A}_0$, $\{e_i\}$ is a basis for $\mathcal{A}_1$, and $\{T_i\}$ is a basis for $\mathcal{A}_2$.
\end{prop}

\begin{rmk}
     More generally, the Eagon-Northcott complex resolves a determinantal ring defined by the $k \times k$ minors of an $n \times m$ generic matrix over the ambient polynomial ring \cite{cawavtag}.  In our case, the Hilbert-Burch and Eagon-Northcott complexes coincide.  The Eagon-Northcott complex admits the structure of a dg algebra, computed by Srinivasan in \cite{hema}.  
\end{rmk}

\begin{defn}\label{Goloddef}
    Given a surjective homomorphism of local rings $(S,\m_S,k_S) \to (S/J,\m_{S/J},k_S)$ and a finitely-generated $S/J$-module $N$, there is a coefficientwise inequality of Poincar\'e series
    \[
    P_N^{S/J}(t) \preceq \dfrac{P^S_N(t)}{1-t(P^S_{S/J}(t)-1)}.
    \]
    In the case that equality holds, $N$ is called \textit{Golod}.  We call the ring $S/J$ \textit{Golod} if $k_S$ is Golod.  
\end{defn}

When the minimal free resolution of $S/J$ over $S$ is a dg algebra, the next proposition provides an equivalent condition for $S/J$ to be Golod.  We phrase it in terms of $R$, $Q$, and $\mathcal{A}$.  

\begin{prop}[Proposition 5.2.4, \cite{IFR}]\label{golodcon}
The ring $R$ is Golod if and only if $\mathcal{A}_{\geq 1}\cdot\mathcal{A}_{\geq 1} \subseteq\m_Q\mathcal{A}$, where $\m_Q$ is the (homogeneous) maximal ideal of $Q$.  Similarly, an $R$-module $N$ is Golod if and only if its minimal $Q$-free resolution $\mathcal{U}$ is a dg module over $\mathcal A$ with structure satisfying that the action of $\mathcal{A}_{\geq 1}$ on $\mathcal{U}$ lands in $\m_Q\mathcal{U}$.   
\end{prop}

\begin{cor}\label{Rgol}
    The ring $R$ is Golod.  
\end{cor}

\begin{proof}
    The dg algebra structure of the Hilbert-Burch complex satisfies $\mathcal{A}_{\geq 1}\cdot\mathcal{A}_{\geq 1} \subseteq\m_Q\mathcal{A}$.  
\end{proof}

The construction of the relative bar resolution via dg structures due to Iyengar \cite{Iyengar} is the main tool in this work.  Given a dg algebra $A$ over a ring $S$ resolving $S/J$ and a dg $A$-module $U$ over $S$ resolving an $S/J$-module $N$, set
\[
B_r=\bigoplus_{i_1+\ldots+i_p+j+p=r} S/J \otimes_S \overline{A_{i_1}} \otimes_S \ldots \otimes_S \overline{A_{i_p}} \otimes_S U_j, 
\]
where $\overline{A}$ is the complex of $S$-modules $\coker(S \to A)$.  In our setting, where $A_0=S$, this says that $\overline{A}=\Sigma A_{\geq 1}$, the shift of the positive part of $A$ by one homological degree. 

\begin{defn}
The \textit{relative bar resolution} has $B_r$ as its $r$th term and differential
\begin{align*}
\partial(s\otimes a_1 \otimes \ldots \otimes a_p \otimes u)=&\pm\sum_{k=1}^p s \otimes a_1 \otimes \ldots \otimes \partial^A(a_k) \ldots \otimes a_p\otimes u\\
&\pm s \otimes a_1 \otimes \ldots \otimes a_p \otimes \partial^{U}(u)\\
&\pm\sum_{k=1}^{p-1} s \otimes a_1 \otimes \ldots \otimes a_k\cdot a_{k+1} \otimes \ldots \otimes a_p\otimes u\\
&\pm s \otimes a_1 \otimes \ldots \otimes a_p \cdot u,
\end{align*}
where $s \in S/J$, $a_\ell \in \overline{A_{i_\ell}}$, and $u \in U_j$.  Exact signs are suppressed, and they are unimportant to the moral of the story.
\end{defn}

\begin{rmk}
One can define the relative bar resolution more generally using $A_\infty$-structures (see, for example, \cite{Burke}).  For our results, we need only the definition in the dg algebra setting. 
\end{rmk}

\begin{rmk}
    The relative bar resolution is a resolution of $N$ over $S/J$.  By considering the construction of the differentials in the relative bar resolution via dg structures and Proposition \ref{golodcon}, one observes that if $S/J$ is Golod, $N$ is an $S/J$-module which is Golod, and both $A$ and $U$ are minimal, then the relative bar resolution resolving $N$ over $S/J$ is minimal.  Thus, in this setting, the Poincar\'e series for $N$ over $S/J$ may be realized as in Definition \ref{Goloddef} when equality holds.  We will show that this is the case in our setting, where $S=Q$, $S/J=R$, and $N$ is the cokernel of the Jacobian matrix $J^T_R$ whose kernel is $\Der_{R \mid k}$ (to be defined in a moment).  The following will be expanded upon in Corollary \ref{poin}:
    \[
    P^{R}_N(t)=\dfrac{P^{Q}_N(t)}{1-t(P^{Q}_{R}(t)-1)}.
    \]
\end{rmk}

The module we will be resolving with a truncation of the relative bar resolution is the module of $k$-linear derivations on $R$.

\begin{defn}
    A $k$-linear \textit{derivation} on a $k$-algebra $S$ is a $k$-linear function $\delta:S \to S$ satisfying the \textit{Leibniz rule}:
    \[
    \delta(ab)=\delta(a)b+a\delta(b),
    \]
    for all $a,b \in S$.  The collection of $k$-linear derivations on $S$ forms an $S$-module under addition of functions and $S$-scalar multiplication.  Denote this $S$-module $\Der_{S \mid k}$.  
\end{defn}

\begin{ex}\label{eu}
    A $k$-linear derivation $\delta$ on the polynomial $k$-algebra $k[x_1,\ldots,x_s]$ can be written
    \[
    \delta=\sum_{j=1}^m\alpha_j \dfrac{\partial}{\partial x_{j}},
    \]
where $\alpha_j \in k[x_1,\ldots,x_s]$ and $\dfrac{\partial}{\partial x_{j}}$ is partial differentiation with respect to $x_j$.  The prototypical example of a derivation is the Euler derivation
\[
E:=\sum_{j=1}^s x_j\dfrac{\partial}{\partial x_j},
\]
which satisfies $E(f):=\deg(f)\cdot f$ for homogeneous $f$.
\end{ex}

Write $J^T_Q$ for the transpose of the Jacobian matrix of $I_n$ in $Q$.  That is, the $i$th row of $J^T_Q$ is indexed by $F_i$ and the columns of $J^T_Q$ are indexed by the indeterminates $x_{i,j}$ according to the rule $x_{i,j} \prec x_{i',j'}$ precisely when $i \leq i'$ and if $i=i'$, $j \leq j'$, i.e.,  
\[
x_{1,1} \prec x_{1,2} \prec \ldots \prec x_{1,n+1} \prec x_{2,1} \prec x_{2,2} \prec \ldots \prec x_{n,n+1},
\]
so that the $(i,(j,k))$-entry of $J^T_Q$ is $\dfrac{\partial F_i}{\partial x_{j,k}}$.  We define $J^T_R$ to be $J^T_Q \otimes_Q R$ and keep the above labeling convention throughout this work.  

\begin{ex}
In the $2 \times 3$ case,
\[
X=\begin{bmatrix}
    x_{1,1}&x_{1,2}&x_{1,3}\\
    x_{2,1}&x_{2,2}&x_{2,3}
\end{bmatrix},
\]
$F_1=x_{1,2}x_{2,3}-x_{1,3}x_{2,2}$, $F_2=x_{1,1}x_{2,3}-x_{1,3}x_{2,1}$, and $F_3=x_{1,1}x_{2,2}-x_{1,2}x_{2,1}$, so that
\[
J^T_Q=\begin{bmatrix}
    0&x_{2,3}&-x_{2,2}&0&-x_{1,3}&x_{1,2}\\
    x_{2,3}&0&-x_{2,1}&-x_{1,3}&0&x_{1,1}\\
    x_{2,2}&-x_{2,1}&0&-x_{1,2}&x_{1,1}&0
\end{bmatrix}.
\]
\end{ex}

\begin{rmk}\label{trivial}
    Notice that $\dfrac{\partial F_i}{\partial x_{k,i}}=0$ for all $i$ and $k$, since to compute $F_i$, the $i$th column of $X$ is deleted.  
\end{rmk}

\section{Useful Lemmas}\label{lemmas}

In this section we further explore the relationship between minors and partial derivatives, developing some useful identities.  These identities are involved primarily in the proofs of Lemmas \ref{cx} and \ref{JT}, which combine with the Buchsbaum-Eisenbud Acyclicity Theorem (Theorem \ref{BEA}) to produce the minimal $Q$-free resolution $\mathcal{U}$ of $\coker J^T_R$, and in the proof of the dg action of the Hilbert-Burch dg algebra $\mathcal{A}$ (Proposition \ref{EN}) on $\mathcal{U}$ (Lemma \ref{act}).  

\begin{rmk}\label{hilbertmult} Recall that for lists of positive integers $\alpha$ and $\beta$, we write $X_\alpha^\beta$ for the matrix obtained by deleting the columns with labels in $\alpha$ and rows with labels in $\beta$.  To compute $F_i$, one starts with $X$ and deletes column $i$ to obtain the matrix $X_i$.  The determinant of $X_i$ is $F_i$, and one can compute this determinant via Laplace expansion along the column of $X_i$ labeled $j$ (when $i \neq j$).  Finding $\dfrac{\partial F_i}{\partial x_{k,j}}$ thus amounts to finding the determinants of $X_{i,j}^k$ for varying $k$, i.e., the polynomial coefficient of $x_{k,j}$ in $F_i$.  However, if $i<j$, the sign associated to $x_{k,j}$ in $X$, which is $(-1)^{j+k}$, switches to $(-1)^{j+k+1}$ in $X_i$.  In light of this observation, we note that when $i \neq j$, we can rewrite the multiplication in the dg algebra structure of the Hilbert-Burch complex (Proposition \ref{EN}) as
    \[
    e_i \cdot e_j=\begin{cases}
        \sum_{k=1}^n(-1)^{i+j+k}(-1)^{j+k+1}\dfrac{\partial F_i}{\partial x_{k,j}}T_k= \sum_{k=1}^n(-1)^{i+1}\dfrac{\partial F_i}{\partial x_{k,j}}T_k, & \text{if } i<j,\\\\
        \sum_{k=1}^n(-1)^{i+j+k}(-1)^{j+k}\dfrac{\partial F_i}{\partial x_{k,j}}T_k= \sum_{k=1}^n(-1)^{i}\dfrac{\partial F_i}{\partial x_{k,j}}T_k, & \text{if } i>j.
        \end{cases} 
    \]  
\end{rmk}

\begin{rmk}\label{swap}
    We claim that
    \[
    \dfrac{\partial F_i}{\partial x_{k,j}}=(-1)^{i+j+1}\dfrac{\partial F_j}{\partial x_{k,i}}.
    \]
    Without loss of generality, let $i>j$.  We have
    \begin{align*}
        (-1)^{k+j}\dfrac{\partial F_i}{\partial x_{k,j}}=\det X_{i,j}^k=(-1)^{k+i+1}\dfrac{\partial F_j}{\partial x_{k,i}},
    \end{align*}
    from which the claim follows.  
\end{rmk}


\begin{defn}
    The \textit{rank} of a matrix $M$ (or of a map represented by a matrix $M$) is the size of the largest non-vanishing minor of $M$.  
\end{defn}

\begin{lem}\label{rankJ^T}
The rank of $J^T_Q$ is $n+1$.  
\end{lem}

\begin{proof}

Consider the maximal minor of $J^T_Q$ determined by the last column and the $n$ columns with first entries $\dfrac{\partial F_1}{\partial x_{u,1}}$, $u \in \{1,\ldots,n\}$.  These entries are zero, and so the minor in question is a product of $\pm\dfrac{\partial F_1}{\partial x_{n,n+1}}$ and the determinant of a matrix $Y$ with entries $\dfrac{\partial F_{r \neq 1}}{\partial x_{u,1}}$, $u \in \{1,\ldots,n\}$.

Specialize by sending $x_{i,i+1}$ to itself for $i=1,\ldots,n-1$ and $x_{i,j}$ to zero for $i\neq j-1$.  The matrix $Y$ becomes diagonal under this specialization, since the entries that were $\dfrac{\partial F_{u+1}}{\partial x_{u,1}}$ remain nonzero (since one summand is $\pm \displaystyle\sum_{i\neq u} x_{i,i+1}$) and the entries that were $\dfrac{\partial F_r}{\partial x_{u,1}}$ become zero for all $u$ and $r \neq u+1$.  The entry $\dfrac{\partial F_1}{\partial x_{n,n+1}}$ of $J^T_Q$ remains nonzero under this specialization, as well.  

The determinant of $Y$ under the specialization is thus a product of nonzero entries in a domain, and so the minor we considered at the start of the proof is nonzero.  Thus, $J^T_Q$ has full rank $n+1$. 
\end{proof}

\begin{ex}
    In the $2 \times 3$ case, i.e., the case when $n=2$, we have that $\dfrac{\partial F_1}{\partial x_{2,3}}=x_{1,2}$.  The matrix $Y$ from the proof of Lemma \ref{rankJ^T} is
    \[
    Y=\begin{bmatrix}
        x_{2,3}&-x_{1,3}\\
        x_{2,2}&-x_{1,2}
    \end{bmatrix},
    \]
    and after specialization, we get
    \[
    \begin{bmatrix}
        x_{2,3}&0\\
        0&-x_{1,2}
    \end{bmatrix}.
    \]
\end{ex}

\begin{defn}
The \textit{grade} against $Q$ of a $Q$-ideal $I$ is the length of the longest $Q$-sequence in $I$, denoted $\grade(I,Q)$.  
\end{defn}

\begin{rmk}\label{rmk}
    The grade against $Q$ of the $Q$-ideal generated by the rank-size minors of $J^T_Q$ is at least 1, since this ideal is nonzero by Lemma \ref{rankJ^T} and $Q$ is a domain. 
\end{rmk}

The next lemma provides some identities involving partials and minors that will be used in developing the minimal $Q$-free resolution of $\coker J^T_R$ in Section \ref{rescok} and the differential graded structure appearing in Section \ref{dgastructure}.   

\begin{lem}\label{syzarg}
There are equalities
\begin{align*}
&\sum_{u=1}^{n+1} x_{i,u}\dfrac{\partial F_r}{\partial x_{s,u}}=\begin{cases}
    0, & i \neq s,\\
    F_r, & i=s,
\end{cases} \tag{1}\\\\
    &\sum_{k=1}^n x_{k,i}\dfrac{\partial F_i}{\partial x_{k,j}}=(-1)^{i+j+1}F_j\, \text{ when } i \neq j, \tag{2}\\\\
&\sum_{r \neq i} (-1)^r\dfrac{\partial F_r}{\partial x_{j,i}} x_{\ell,r}=\begin{cases}0, & j \neq \ell,\\(-1)^{i+1}F_i,& j=\ell,\end{cases} \tag{3}\\\\
&\sum_{\ell=1}^n\dfrac{\partial F_i}{\partial x_{\ell,k}}x_{\ell,t}=0 \text{ when $i$, $k$, and $t$ are distinct}.  \tag{4}  
\end{align*}

\end{lem}

\begin{proof}

    To see the case of (1) when $i \neq s$, let $Y$ be the matrix obtained from $X_r$ by replacing the row indexed $s$ with a second row indexed $i$.  The determinant of $Y$ is zero, and one way to compute $\det Y$ is via Laplace expansion along one of the rows indexed by $i$ in $Y$.  Explicitly, this computation is
    \[
    0=\det Y=\pm\sum_{u=1}^{n+1} x_{i,u}\dfrac{\partial F_r}{\partial x_{s,u}}.
    \]
    This yields the case of (1) in which $i \neq s$. 

    When $i=s$, the result follows from the observation that the sum in (1) is the calculation of $F_r$ via Laplace expansion along the row of $X_r$ indexed $s$.   

    For equality (2), we examine $x_{k,i}\dfrac{\partial F_i}{\partial x_{k,j}}$:
    \[
    x_{k,i}\dfrac{\partial F_i}{\partial x_{k,j}}=\begin{cases}
        (-1)^{j+k+1}x_{k,i}\det X_{i,j}^k, & i<j,\\
        (-1)^{j+k}x_{k,i}\det X_{i,j}^k, & i>j.
    \end{cases}
    \]
    By Laplace expansion to compute $F_j$ along the column of $X_j$ indexed $i$,
    \[
    \sum_{k=1}^n(-1)^{i+k+1}x_{k,i}\det X_{i,j}^k=F_j
    \]
    when $i>j$, and
    \[
    \sum_{k=1}^n(-1)^{i+k}x_{k,i}\det X_{i,j}^k=F_j
    \]
    when $i<j$.  Equality (2) follows:
     \[
    \sum_{k=1}^n x_{k,i}\dfrac{\partial F_i}{\partial x_{k,j}}=(-1)^{i+j+1}F_j.
    \]

    Alternatively, equation (2) follows from Remark \ref{swap} and Laplace expansion as in the $i=s$ case of equation (1).

    To see the case of (3) when $j \neq \ell$, let $Y$ be the matrix obtained from $X_i$ by replacing the row indexed by $j$ with a second row indexed by $\ell$.  By Laplace expansion along the row in $Y$ originally indexed by $\ell$, we get
    \begin{align*}
    0=\det Y&=\sum_{r=1}^{i-1} (-1)^{\ell+r}x_{\ell,r}\det X_{i,r}^j+\sum_{r=i+1}^n (-1)^{\ell+r+1}x_{\ell,r}\det X_{i,r}^j.
    \end{align*}

    When $i>r$,
    \[
    \dfrac{\partial F_r}{\partial x_{j,i}}=(-1)^{i+j+1} \det X_{i,r}^j,
    \]
    and when $i<r$, 
    \[
    \dfrac{\partial F_r}{\partial x_{j,i}}=(-1)^{i+j} \det X_{i,r}^j.
    \]
    So,
    \begin{align*}
        \det Y&=\sum_{r=1}^{i-1} (-1)^{\ell+r+i+j+1}x_{\ell,r}\dfrac{\partial F_r}{\partial x_{j,i}}+\sum_{r=i+1}^n (-1)^{\ell+r+i+j+1}x_{\ell,r}\dfrac{\partial F_r}{\partial x_{j,i}}\\
        &=(-1)^{\ell+i+j+1}\sum_{r \neq i}(-1)^rx_{\ell,r}\dfrac{\partial F_r}{\partial x_{j,i}}.
    \end{align*}
    Thus,
    \[
    \sum_{r\neq i}\left((-1)^{r}\dfrac{\partial F_r}{\partial x_{j,i}}x_{\ell,r}b_{\ell,i}\right)=0.  
    \]

    To see the case of (3) when $j=\ell$, consider the determinant of $X_i$ as computed via Laplace expansion along the row indexed by $j$:
    \begin{align*}
    F_i=\det X_i&=\sum_{r=1}^{i-1} (-1)^{j+r}x_{j,r}\det X_{i,r}^j+\sum_{r=i+1}^n (-1)^{j+r+1}x_{j,r}\det X_{i,r}^j.
    \end{align*}
    
    When $i>r$,
    \[
    \dfrac{\partial F_r}{\partial x_{j,i}}=(-1)^{i+j+1} \det X_{i,r}^j,
    \]
    and when $i<r$, 
    \[
    \dfrac{\partial F_r}{\partial x_{j,i}}=(-1)^{i+j} \det X_{i,r}^j.
    \]
    Thus,
    \begin{align*}
    F_i&=\sum_{r=1}^{i-1} (-1)^{i+r+1}x_{j,r}\dfrac{\partial F_r}{\partial x_{j,i}}+\sum_{r=i+1}^n (-1)^{i+r+1}x_{j,r}\dfrac{\partial F_r}{\partial x_{j,i}}\\
    &=(-1)^{i+1}\sum_{r\neq i}(-1)^{r}\dfrac{\partial F_r}{\partial x_{j,i}}x_{j,r},
    \end{align*}
    from which (3) follows.  

    To see (4), let $Y$ be the matrix obtained from $X_i$ by replacing the column indexed $k$ by a second column indexed $t$.  By Laplace expansion along one of the columns indexed $t$ in $Y$,
    \[
    0=\det Y=\pm\sum_{\ell=1}^n\dfrac{\partial F_i}{\partial x_{\ell,k}}x_{\ell,t}.
    \]
    
    \end{proof}

\begin{ex}
In the $2 \times 3$ case with $i=2$ and $j=1$, equation (3) in Lemma \ref{syzarg} reads
\begin{align*}
\sum_{r \neq 2}(-1)^r\dfrac{\partial F_r}{\partial x_{1,2}}x_{1,r}&=-\dfrac{\partial F_1}{\partial x_{1,2}}x_{1,1}-\dfrac{\partial F_3}{\partial x_{1,2}}x_{1,3}\\
&=-x_{2,3}x_{1,1}+x_{2,1}x_{1,3}\\
&=-F_2.
\end{align*}
\end{ex}

\section{Resolving the Cokernel over $Q$}\label{rescok}

In the introduction, we noted that $\Der_{R \mid k}$ is isomorphic to the kernel of $J^T_R$, i.e., there is an exact sequence of $R$-modules
\[
0 \to \Der_{R \mid k} \to R^{n(n+1)} \xrightarrow{J^T_R} R^{n+1} \to \coker J^T_R \to 0.  
\]
Resolving $\Der_{R \mid k}$ over $R$ thus amounts to resolving $\coker J^T_R$ over $R$.  To do so, we first resolve $\coker J^T_R$ over $Q$, and that is the content of this section. 
 Using the resolution we develop in this section (see Theorem \ref{cokresn}) together with the Hilbert-Burch complex from Proposition \ref{EN}, we construct the relative bar resolution of $\coker J^T_R$ over $R$ in Section \ref{resder}.  ``Chopping off" the first term of this resolution results in the minimal $R$-free resolution of $\Der_{R \mid k}$.    

We start by examining what will be the second differential in the $Q$-free resolution of $\coker J^T_R$.  Let $A_k$ be the $(n+1) \times (n-1)$ matrix obtained from the $(n+1) \times n$ matrix
    \[
    X^T=\begin{bmatrix}
        x_{1,1} & x_{2,1} & \ldots & x_{n,1}\\
        x_{1,2} & x_{2,2} & \ldots & x_{n,2}\\
        \vdots & \vdots & \vdots & \vdots\\
        x_{1,n+1} & x_{2,n+1} & \ldots & x_{n,n+1}\\
    \end{bmatrix}
    \] 
    by removing the $k$th column.  Let $\vec{e}_\ell$ be the $\ell$th $n(n+1) \times 1$ standard basis vector and let $B_\ell$ be the $n(n+1) \times 1$ column matrix 
    \[
    \sum_{s=1}^{n+1}\left( -x_{1,s}\vec{e}_s+x_{\ell,s}\vec{e}_{(\ell-1)(n+1)+s}\right). 
    \]

   Let $\partial_2:Q^{(n-1)(n+1)} \to Q^{n(n+1)}$ be the map of free $Q$-modules be given by the block matrix
   \[
   \begin{bmatrix}
       A_1&0&\ldots&0&0&\vert&\vert&\vert&\vert\\
       0&A_2&\ldots&0&0&\vert&\vert&\vert&\vert\\
       \vdots&\vdots&\ddots&\vdots&\vdots&B_2&B_3&\ldots&B_n\\
       0&0&\ldots&A_{n-1}&0&\vert&\vert&\vert&\vert\\
       0&0&\ldots&0&A_n&\vert&\vert&\vert&\vert
   \end{bmatrix}
   \]
   
Let $I$ be the ideal of $(n-1)(n+1)\times(n-1)(n+1)$ (i.e., maximal) minors of $\partial_2$.  

\begin{lem}\label{partial2hasgrade2}
   We have that $\grade(I,Q)$ is at least 2. 
\end{lem}

\begin{proof}
     For a matrix $Y$ with entries in $Q$, denote by $I(Y)$ the $Q$-ideal generated by the maximal minors of $Y$.  

    Take an $(n-1)(n+1)$ minor of $\partial_2$ by choosing $n-1$ rows from the collection of $n+1$ rows that intersects $A_1$ and $n$ rows from each of the $n-1$ collections of $n+1$ rows that intersect each $A_k$, $k\geq 2$.  Using column operations, one sees that the matrix corresponding to this minor can be calculated along a block diagonal as an element of $I
    (X^1)\cdot I(X)^{n-1}$.  The collection of minors taken this way generates the ideal $I(X^1)\cdot I(X)^{n-1}$, so that $I(X^1)\cdot I(X)^{n-1} \subseteq I$.  Using that $I(X)^{n-1}\subseteq I(X) \subseteq I(X^1)$, we have
    \[
    \sqrt{I(X^1)\cdot I(X)^{n-1}}=\sqrt{I(X^1)\cap I(X)^{n-1}}=\sqrt{I(X)^{n-1}}=\sqrt{I(X)}=I(X).
    \]  
    We use above that $I(X)$ is prime and thus radical.  We note that $\height \left(\sqrt{I(X^1)\cdot I(X)^{n-1}}\right)=\height \left(I(X^1)\cdot I(X)^{n-1}\right)$ and that $\grade(I(X),Q)=2$ (see, for example, \cite{bh}, Theorem 7.3.1).  Since height and grade coincide in the case of a Cohen-Macaulay ring, $\grade(I,Q) \geq 2$.  
\end{proof}

\begin{ex}
    In the $3 \times 4$ case, the matrix $\partial_2$ is
    \[
    \begin{bmatrix}
        \textcolor{red}{x_{2,1}}&\textcolor{red}{x_{3,1}}&\textcolor{red}{0}&\textcolor{red}{0}&\textcolor{red}{0}&\textcolor{red}{0}&\textcolor{red}{-x_{1,1}}&\textcolor{red}{-x_{1,1}}\\
        x_{2,2}&x_{3,2}&0&0&0&0&-x_{1,2}&-x_{1,2}\\
        \textcolor{red}{x_{2,3}}&\textcolor{red}{x_{3,3}}&\textcolor{red}{0}&\textcolor{red}{0}&\textcolor{red}{0}&\textcolor{red}{0}&\textcolor{red}{-x_{1,3}}&\textcolor{red}{-x_{1,3}}\\
        x_{2,4}&x_{3,4}&0&0&0&0&-x_{1,4}&-x_{1,4}\\
        \textcolor{red}{0}&\textcolor{red}{0}&\textcolor{red}{x_{1,1}}&\textcolor{red}{x_{3,1}}&\textcolor{red}{0}&\textcolor{red}{0}&\textcolor{red}{x_{2,1}}&\textcolor{red}{0}\\
        \textcolor{red}{0}&\textcolor{red}{0}&\textcolor{red}{x_{1,2}}&\textcolor{red}{x_{3,2}}&\textcolor{red}{0}&\textcolor{red}{0}&\textcolor{red}{x_{2,2}}&\textcolor{red}{0}\\
        \textcolor{red}{0}&\textcolor{red}{0}&\textcolor{red}{x_{1,3}}&\textcolor{red}{x_{3,3}}&\textcolor{red}{0}&\textcolor{red}{0}&\textcolor{red}{x_{2,3}}&\textcolor{red}{0}\\
        0&0&x_{1,4}&x_{3,4}&0&0&x_{2,4}&0\\
        \textcolor{red}{0}&\textcolor{red}{0}&\textcolor{red}{0}&\textcolor{red}{0}&\textcolor{red}{x_{1,1}}&\textcolor{red}{x_{2,1}}&\textcolor{red}{0}&\textcolor{red}{x_{3,1}}\\
        0&0&0&0&x_{1,2}&x_{2,2}&0&x_{3,2}\\
        \textcolor{red}{0}&\textcolor{red}{0}&\textcolor{red}{0}&\textcolor{red}{0}&\textcolor{red}{x_{1,3}}&\textcolor{red}{x_{2,3}}&\textcolor{red}{0}&\textcolor{red}{x_{3,3}}\\
        \textcolor{red}{0}&\textcolor{red}{0}&\textcolor{red}{0}&\textcolor{red}{0}&\textcolor{red}{x_{1,4}}&\textcolor{red}{x_{2,4}}&\textcolor{red}{0}&\textcolor{red}{x_{3,4}}\\
    \end{bmatrix},
    \]
    in which we have highlighted a maximal minor described in the proof of Lemma \ref{partial2hasgrade2}.  
\end{ex}

\begin{rmk}\label{rankdel2}
In Lemma \ref{partial2hasgrade2}, the minors we take are products of nonzero elements of $Q$, which is a domain (the factors are determinants of matrices of distinct indeterminates, hence nonzero).  Thus, the map $\partial_2$ has full rank $(n-1)(n+1)$.  
\end{rmk}

Having said something about the map we will claim as the second differential in the resolution of $\coker J^T_R$ over $Q$, we turn now to the rest of the resolution.  We begin with a lemma that will aid in showing that the resolution is indeed a complex.  

\begin{lem}\label{lem}
Let $S=k[x_1,\ldots,x_m]$ and $f\in S$.  Let $\delta \in \Der_{S \mid k}$ have coordinate vector $\vec{\alpha}$, and let $J^T_S(f)$ be the transpose of the Jacobian of $f$ over $S$.  The composition $J^T_S(f) \circ \vec{\alpha}=\vec{0}$ over $S$ if and only if $\delta(f)=0$.  
\end{lem}

\begin{proof}
    Say $\delta=\sum_{j=1}^m\alpha_j\dfrac{\partial}{\partial x_j}$.  We have that $J^T_S(f) \circ \vec{\alpha}$ is equal to
    \[
    \sum_{j=1}^m\dfrac{\partial f}{\partial x_j}\alpha_j=\sum_{j=1}^m\alpha_j\dfrac{\partial f}{\partial x_j}=\delta(f).
    \]
\end{proof}

\begin{lem}\label{cx}
    The sequence of maps and free modules
    \[
    Q^{(n-1)(n+1)} \xrightarrow{\partial_2} Q^{n(n+1)} \xrightarrow{J^T_Q} Q^{n+1}
    \]
    is a complex.
\end{lem}

\begin{proof}
We want to show that $J^T_Q \circ \partial_2=0$.  By Lemma \ref{lem}, it suffices to show that the columns of $\partial_2$ correspond to $k$-linear derivations on $Q$ that send the polynomials $F_r$ (the generators of $I_n)$ to zero.  Define
    \[
    V_{i,s}:=\sum_{u=1}^{n+1}x_{i,u}\dfrac{\partial}{\partial x_{s,u}}.
    \]
    For $i \neq s$, $V_{i,s}$ is the derivation corresponding to the column of $\partial_2$ having first nonzero entry $x_{i,1}$.  We will use $B_\ell$ to denote the derivation corresponding to the column $B_\ell$.  Notice that $B_\ell=-V_{1,1}+V_{\ell,\ell}$.  

    That $V_{i,s}(F_r)=0$ when $i \neq s$ and $V_{s,s}(F_r)=F_r$ is the content of Lemma \ref{syzarg}, equation (1).  Since $V_{s,s}(F_r)=F_r$, we get that $B_\ell(F_r)=0$.      

    By Lemma \ref{lem}, $J^T_Q \circ \partial_2=0$. 
\end{proof}

We could extend the complex from Lemma \ref{cx} to the right to include a surjection onto $\coker J^T_Q$.  However, we would like to resolve $\coker J^T_R$ over $Q$.  The next lemma states that this is not an issue.    

\begin{lem}\label{JT}
    As $Q$-modules, $\coker J^T_Q \cong \coker J^T_R$. 
\end{lem}

\begin{proof}
    We have that
    \[
    \coker J^T_R = \frac{R^{n+1}}{\im J^T_R} \cong \frac{(Q/I)^{n+1}}{\im J^T_R} \cong \frac{Q^{n+1}}{\la F_i\vec e_j \ra+\im J^T_Q},
    \]
    where $\{\vec e_j\}$ is a basis for the $j$th copy of $Q$ and $\la F_i\vec e_j \ra$ is the ideal of $Q^{n+1}$ generated by the $F_i\vec e_j$, $i,j=1,\ldots,n+1$.  The module $\dfrac{Q^{n+1}}{\la F_i\vec e_j \ra+\im J^T_Q}$ is isomorphic to the cokernel of a matrix equal to $J^T_Q$ with adjoined columns that are standard basis vectors scaled by the $F_i$.  Call this matrix $D$.  The lemma follows if we can show that the columns of $D$ are $Q$-linear combinations of the columns of $J^T_Q$.  
    
    Write $b_{i,j}$ for the column of $J^T_Q$ corresponding to the partial derivative with respect to $x_{i,j}$.  When $i \neq j$, we claim that
    \[
    \sum_{k=1}^n x_{k,i}b_{k,j}=(-1)^{i+j+1}F_j\vec e_i.
    \]
    The $r$th entry of the vector $\sum_{k=1}^n x_{k,i}b_{k,j}$ is
    \[
    \sum_{k=1}^nx_{k,i}\dfrac{\partial F_r}{\partial x_{k,j}}=\begin{cases}
        (-1)^{i+j+1}F_j, & i=r,\\
        0, & i \neq r.
    \end{cases}
    \]
    The $i=r$ case is Lemma \ref{syzarg}, equation (2).  The $i \neq r$ case is Lemma \ref{syzarg}, equation (4).  

    When $i=j$, we claim that
    \begin{align*}
    \sum_{u \neq i}x_{i,u}b_{i,u}-\sum_{k \neq i}x_{k,i}b_{k,i}=F_i\vec e_i.  
    \end{align*}
    The $r$th entry of $\sum_{u \neq i}x_{i,u}b_{i,u}-\sum_{k\neq i}x_{k,i}b_{k,i}$ is
    \begin{align*}
    \sum_{u \neq i}x_{i,u}\dfrac{\partial F_r}{\partial x_{i,u}}-\sum_{k\neq i}x_{k,i}\dfrac{\partial F_r}{\partial x_{k,i}}.
    \end{align*}
    When $r=i$, this sum is $F_i$ by Laplace expansion computing $F_i$ along the row of $X_i$ indexed by $i$.  When $r \neq i$, 
    \begin{align*}
    \sum_{u \neq i}x_{i,u}\dfrac{\partial F_r}{\partial x_{i,u}}-\sum_{k\neq i}x_{k,i}\dfrac{\partial F_r}{\partial x_{k,i}}&=\sum_{u=1}^{n+1}x_{i,u}\dfrac{\partial F_r}{\partial x_{i,u}}-x_{i,i}\dfrac{\partial F_r}{\partial x_{i,i}}-\left(\sum_{k=1}^nx_{k,i}\dfrac{\partial F_r}{\partial x_{k,i}}-x_{i,i}\dfrac{\partial F_r}{\partial x_{i,i}}\right)\\
    &=F_r-x_{i,i}\dfrac{\partial F_r}{\partial x_{i,i}}-F_r+x_{i,i}\dfrac{\partial F_r}{\partial x_{i,i}}\\
    &=0
    \end{align*}
\end{proof}

To show that what we claim is the $Q$-free resolution of $\coker J^T_R$ is in fact acyclic, we call upon the Buchsbaum-Eisenbud Acyclicity Theorem.  

\begin{thm}[Buchsbaum-Eisenbud Acyclicity \cite{BEA}]\label{BEA}
    Let 
    \[
    0 \longrightarrow F_s \xrightarrow{d_s} F_{s-1} \longrightarrow \ldots \longrightarrow F_1 \xrightarrow{d_1} F_0 \tag{$\mathcal{F}$}
    \]
    be a complex of free $Q$-modules in which $F_i=Q^{f_i}$.  Let $r_i$ be the rank of the matrix representing $d_i$.  By convention, $r_{s+1}=0$.  Then, $\mathcal{F}$ is acyclic if and only if the following two conditions hold
    \begin{itemize}
        \item[(i)] $f_i=r_i+r_{i+1}$ for all $1 \leq i \leq s$,
        \item[(ii)] for all $1 \leq i \leq s$, we have that $\grade(I(d_i),Q) \geq i$, where $I(d_i)$ is the ideal of $Q$ generated by to $r_i \times r_i$ minors of $d_i$.  
    \end{itemize}
\end{thm}

To end this section, we bring its contents together to obtain the minimal $Q$-free resolution of $\coker J^T_R$.  

\begin{thm}\label{mfrcokerjt}
    The minimal $Q$-free resolution of $\coker J^T_R$ is
    \[
    0 \longrightarrow Q^{(n-1)(n+1)} \xrightarrow{\partial_2} Q^{n(n+1)} \xrightarrow{J^T_Q} Q^{n+1} \longrightarrow \coker J^T_R \longrightarrow 0 \tag{$\mathcal{U}$}
    \]
\end{thm}

\begin{proof}
    By Lemma \ref{cx} and Lemma \ref{JT}, $\mathcal{U}$ is a complex.  By Lemma \ref{rankJ^T} and Remark \ref{rankdel2}, the ranks of $J^T$ and $\partial_2$ are maximal.  By Lemma \ref{partial2hasgrade2} and Remark \ref{rmk}, the grades of the ideals of rank-sized minors of the differentials of $\mathcal{U}$ are large enough so that by the Buchsbaum-Eisenbud Acyclicity Theorem (Theorem \ref{BEA}), we may conclude that $\mathcal{U}$ is a resolution.  The entries of the matrices representing the differentials are in the maximal ideal of $Q$, and so $\mathcal{U}$ is indeed minimal.  

\end{proof}
\section{Differential Graded Structure}\label{dgastructure}

Given a dg algebra $\mathcal{A}$ resolving $R$ over $Q$ and a dg $\mathcal{A}$-module $\mathcal{U}$ resolving an $R$-module $N$ over $Q$, one can produce the relative bar resolution resolving $N$ over $R$ (Theorem 1.2, \cite{Iyengar}).  This section develops the dg module structure required when the module $N$ is $\coker J^T_R$.  

We restate the resolutions $\mathcal{A}$ and $\mathcal{U}$ developed earlier, with the free modules in them labeled by their bases.  Recall that $\mathcal{A}$ is the Hilbert-Burch complex (Proposition \ref{EN}):
\[
    0 \longrightarrow \underset{\{T_i\}}{Q^n} \xrightarrow{\varphi_2} \underset{\{e_i\}}{Q^{n+1}} \xrightarrow{\varphi_1} \underset{\{1\}}{Q} \longrightarrow R \longrightarrow 0 \tag{$\mathcal{A}$},
    \]
    where $\varphi_1=\begin{bmatrix}
        F_1&-F_2&\ldots&(-1)^{i+1}F_i&\ldots& (-1)^{n+2}F_{n+1}
    \end{bmatrix}$ and $\varphi_2=X^T$.  Furthermore, the resolution $\mathcal{A}$ is a dg algebra under the multiplication
    \begin{align*}
    e_i \cdot e_j=-(e_j \cdot e_i)&=\begin{cases}
   \displaystyle\sum_{k=1}^n(-1)^{i+j+k}\det X^k_{i,j}T_k, & \text{if } i < j,\\\\
        0, & \text{if } i=j,
    \end{cases}\\\\
    e_i \cdot T_j&=0,\\\\
    T_i \cdot T_j&=0.
    \end{align*}
   Also recall that $\mathcal{U}$ is the minimal $Q$-free resolution of $\coker J^T_R$ (Theorem \ref{mfrcokerjt}):
    \[
    0 \longrightarrow \underset{\{c_{i,j}\}}{Q^{(n-1)(n+1)}} \xrightarrow{\partial_2} \underset{\{b_{j,k}\}}{Q^{n(n+1)}} \xrightarrow{J^T_Q} \underset{\{a_i\}}{Q^{n+1}} \longrightarrow \coker J^T_R \longrightarrow 0 \tag{$\mathcal{U}$}.
    \]
    Restated in more detail, we have   \\

\hspace{5mm} $\bullet$ $\{1\}$ is a basis for $\mathcal{A}_0$\\

\hspace{5mm} $\bullet$ $\{e_i\}$ is a basis for $\mathcal{A}_1$\\

\hspace{5mm} $\bullet$ $\{T_i\}$ is a basis for $\mathcal{A}_2$\\

\hspace{5mm} $\bullet$ $\{a_i\}$ is a basis for $\mathcal{U}_0$\\

\hspace{5mm} $\bullet$ $\{b_{j,k}\}$, $1 \leq j \leq n$, $1 \leq k \leq n+1$, is a basis for $\mathcal{U}_1$\\

\hspace{5mm} $\bullet$ $\partial(b_{i,j})=\displaystyle\sum_{\ell=1}^{n+1}\dfrac{\partial F_\ell}{\partial x_{i,j}}a_\ell$\\

\hspace{5mm} $\bullet$ $\{c_{i,j}\}$, $1 \leq i,j \leq n$, $(i,j) \neq (1,1)$, is a basis for $\mathcal{U}_2$; by convention, we take $c_{1,1}=0$\\

\hspace{5mm} $\bullet$ $\partial (c_{i,j})=\displaystyle\sum_{\ell=1}^{n+1}x_{i,\ell}b_{j,\ell}-\delta_{i,j}\sum_{k=1}^{n+1}x_{1,k}b_{1,k}$, that is, $c_{i,j}$ corresponds

\hspace{5mm} to the column of $\partial_2$ in which $x_{i,k}$ is in the $j$th block of $n+1$ entries\\

\begin{lem}\label{act}
    The resolution $\mathcal{U}$ admits the structure of a dg module over the dg algebra $\mathcal{A}$.  Such a structure is given by
    \begin{align*}
    e_i \cdot a_j&=\begin{cases}
        (-1)^{j}\left(\displaystyle\sum_{k=1}^nx_{k,j}b_{k,i}\right), & \text{if } i \neq j,\\\\
        (-1)^{i+1}\left(\displaystyle\sum_{k \neq i} x_{1,k}b_{1,k}-\displaystyle\sum_{\ell \neq 1} x_{\ell,i}b_{\ell,i}\right), & \text{if } i=j,
    \end{cases}\\\\
    e_i\cdot b_{j,k}&=(-1)^{i+1}\displaystyle\sum_{\ell=1}^n\dfrac{\partial F_i}{\partial x_{\ell,k}}c_{\ell,j},\\\\
    e_i\cdot c_{j,k}&=0,\\\\
    T_i\cdot a_j&=(-1)^j\sum_{k=1}^nx_{k,j}c_{i,k},\\\\
    T_i\cdot b_{j,k}&=0,\\\\
    T_i\cdot c_{j,k}&=0.\\\\
    \end{align*}
\end{lem}

Before the proof, we make a remark and give an example.  

\begin{rmk}
    Notice that in the action $e_i \cdot a_i$, one could take the sums to be over all $k$ and all $\ell$, respectively, since $x_{1,i}b_{1,i}-x_{1,i}b_{1,i}=0$. Additionally, there is nothing special about the number one here.  We could have instead defined
    \[
    e_i\cdot a_i=(-1)^{i+1}\left(\left(\displaystyle\sum_{k \neq i} x_{s,k}b_{s,k}\right)-\left(\displaystyle\sum_{\ell \neq s} x_{\ell,i}b_{\ell,i}\right)\right)
    \]
    for any fixed $s \in \{1,\ldots,n\}$.  This is because the difference
    \begin{align*}
    &(-1)^{i+1}\left(\left(\displaystyle\sum_{k \neq i} x_{s,k}b_{s,k}\right)-\left(\displaystyle\sum_{\ell \neq s} x_{\ell,i}b_{\ell,i}\right)\right)-(-1)^{i+1}\left(\left(\displaystyle\sum_{k \neq i} x_{1,k}b_{1,k}\right)-\left(\displaystyle\sum_{\ell \neq 1} x_{\ell,i}b_{\ell,i}\right)\right)\\
    &=(-1)^{i+1}\left[\left(\displaystyle\sum_{k \neq i} x_{s,k}b_{s,k}\right)-\left(\displaystyle\sum_{k \neq i} x_{1,k}b_{1,k}\right)-x_{1,i}b_{1,i}+x_{s,i}b_{s,i}\right]\\
    &=(-1)^{i+1}\left[\left(\displaystyle\sum_{k=1}^{n+1} x_{s,k}b_{s,k}\right)-\left(\displaystyle\sum_{k=1}^{n+1} x_{1,k}b_{1,k}\right)\right]\\
    &=\partial((-1)^{i+1}c_{s,s})
    \end{align*}
    is a cycle.  
\end{rmk}

\begin{ex}
    Let $n=3$, $i=1$, and $j=k=2$.  Then,
    \begin{align*}
    e_1 \cdot b_{2,2}&=(-1)^2\sum_{\ell=1}^n\dfrac{\partial F_1}{\partial x_{\ell,2}}c_{\ell,2}\\
    &=\dfrac{\partial F_1}{\partial x_{1,2}}c_{1,2}+\dfrac{\partial F_1}{\partial x_{2,2}}c_{2,2}+\dfrac{\partial F_1}{\partial x_{3,2}}c_{3,2}\\
    &=(x_{2,3}x_{3,4}-x_{2,4}x_{3,3})c_{1,2}-(x_{1,3}x_{3,4}-x_{1,4}x_{3,3})c_{2,2}+(x_{1,3}x_{2,4}-x_{1,4}x_{2,3})c_{3,2}.
    \end{align*}
\end{ex}


\begin{proof}[Proof of Lemma \ref{act}]
    We show a number of computations in the proof.  While similar in flavor, we include as much as we do since different cases require different tricks.  

    The fact that $\mathcal{U}$ is a dg module over $\mathcal{A}$ follows from Proposition 2.2.5 in \cite{IFR}, but the explicit structure does not.  That the products $e_i \cdot c_{j,k}$, $T_i \cdot b_{j,k}$, and $T_i \cdot c_{j,k}$ are all equal to zero comes from the fact that $\mathcal{U}_{\geq 3}=0$.    

    In what follows, an equality is tagged $(n)$ if it follows from Lemma \ref{syzarg}, equation $(n)$. 

    \textbf{Begin $e_ia_j$ action:}
    
    We have that $\partial(e_ia_j)=(-1)^{i+1}F_ia_j$.  When $i \neq j$,
    \begin{align*}
    \partial\left[(-1)^{j}\left(\sum_{k=1}^n x_{k,j}b_{k,i}\right)\right]&=(-1)^{j}\left(\sum_{k=1}^n x_{k,j}\partial(b_{k,i})\right)\\
    &=(-1)^{j}\left(\sum_{k=1}^n x_{k,j}\left(\sum_{\ell=1}^{n+1}\dfrac{\partial F_\ell}{\partial x_{k,i}}a_\ell\right)\right)\\
    &=(-1)^{j}\left(\sum_{\ell=1}^{n+1}\sum_{k=1}^nx_{k,j}\dfrac{\partial F_\ell}{\partial x_{k,i}}a_\ell\right)\\
    &=(-1)^{j}\left(\sum_{k=1}^nx_{k,j}\dfrac{\partial F_j}{\partial x_{k,i}}a_j+\sum_{\ell \neq j}\sum_{k=1}^nx_{k,j}\dfrac{\partial F_\ell}{\partial x_{k,i}}a_\ell\right)\\
    &=(-1)^{j}\left(\sum_{k=1}^nx_{k,j}\dfrac{\partial F_j}{\partial x_{k,i}}a_j\right)\tag{4}\\
    &=(-1)^{j}(-1)^{i+j+1}F_ia_j\tag{2}\\
    &=(-1)^{i+1}F_ia_j.
   \end{align*}

    When $i=j$,
    \begin{align*}
    &\partial\left[(-1)^{i+1}\left(\left(\sum_{k \neq i} x_{1,k}b_{1,k}\right)-\left(\sum_{\ell \neq 1} x_{\ell,i}b_{\ell,i}\right)\right)\right]\\
    &=(-1)^{i+1}\left(\sum_{k \neq i} x_{1,k}\partial(b_{1,k})\right)-\left(\sum_{\ell \neq 1} x_{\ell,i}\partial(b_{\ell,i})\right)\\
    &=(-1)^{i+1}\left(\sum_{k \neq i} x_{1,k}\left(\sum_{u=1}^{n+1}\dfrac{\partial F_u}{\partial x_{1,k}}a_u\right)-\sum_{\ell \neq 1} x_{\ell,i}\left(\sum_{v=1}^{n+1}\dfrac{\partial F_v}{\partial x_{\ell,i}}a_v\right)\right)\\
    &=(-1)^{i+1}\left(\sum_{u=1}^{n+1} \left(F_ua_u-x_{1,i}\frac{\partial F_u}{\partial x_{1,i}}a_u\right)-\sum_{\ell \neq 1} x_{\ell,i}\left(\sum_{v \neq i}\dfrac{\partial F_v}{\partial x_{\ell,i}}a_v\right)\right)\tag{1}\\
    &=(-1)^{i+1}\left(F_ia_i+\sum_{u\neq i} \left(F_ua_u-x_{1,i}\frac{\partial F_u}{\partial x_{1,i}}a_u\right)-\sum_{\ell \neq 1} x_{\ell,i}\left(\sum_{v \neq i}\dfrac{\partial F_v}{\partial x_{\ell,i}}a_v\right)\right)\\
    &=(-1)^{i+1}\left(F_ia_i+\sum_{u\neq i} \left(F_ua_u-x_{1,i}\frac{\partial F_u}{\partial x_{1,i}}a_u\right)-\sum_{v\neq i} \left(F_va_v-x_{1,i}\frac{\partial F_v}{\partial x_{1,i}}a_v\right)\right)\tag{1}\\
    &=(-1)^{i+1}F_ia_i.
    \end{align*}

    \textbf{End $e_ia_j$ action.}

    \textbf{Begin $e_ib_{j,k}$ action:}

    We have that
    \begin{align*}
    \partial(e_ib_{j,k})&=(-1)^{i+1}F_ib_{j,k}-e_i\left(\sum_{r=1}^{n+1}\dfrac{\partial F_r}{\partial x_{j,k}}a_r\right)\\
    &=(-1)^{i+1}F_ib_{j,k}-\sum_{r=1}^{n+1}\dfrac{\partial F_r}{\partial x_{j,k}}e_ia_r\\
    &=(-1)^{i+1}F_ib_{j,k}-\sum_{r\neq i}\dfrac{\partial F_r}{\partial x_{j,k}}e_ia_r-\dfrac{\partial F_i}{\partial x_{j,k}}e_ia_i\\
    &=(-1)^{i+1}F_ib_{j,k}-\sum_{r\neq i}\left(\dfrac{\partial F_r}{\partial x_{j,k}}(-1)^{r}\left(\sum_{\ell=1}^nx_{\ell,r}b_{\ell,i}\right)\right)-\dfrac{\partial F_i}{\partial x_{j,k}}(-1)^{i+1}\left(\sum_{t \neq i} x_{1,t}b_{1,t}-\sum_{u \neq 1} x_{u,i}b_{u,i}\right)
    \end{align*}

    \textbf{Case:} $i=k$.  

    We claim that we may take $e_i \cdot b_{j,i}=(-1)^{i+1}\displaystyle\sum_{\ell=1}^n\dfrac{\partial F_i}{\partial x_{\ell,i}}c_{\ell,j}=0$.  When $i=k$, $\dfrac{\partial F_i}{\partial x_{j,k}}e_ia_i=0$, and so
    \begin{align*}
    \partial(e_ib_{j,i})&=(-1)^{i+1}F_ib_{j,i}-\sum_{r\neq i}\dfrac{\partial F_r}{\partial x_{j,i}}e_ia_r\\
    &=(-1)^{i+1}F_ib_{j,i}-\sum_{r\neq i}\left(\dfrac{\partial F_r}{\partial x_{j,i}}(-1)^{r}\left(\sum_{\ell=1}^nx_{\ell,r}b_{\ell,i}\right)\right)\\
    &=(-1)^{i+1}F_ib_{j,i}-\sum_{r\neq i}\sum_{\ell=1}^n\left((-1)^{r}\dfrac{\partial F_r}{\partial x_{j,i}}x_{\ell,r}b_{\ell,i}\right)\\
    &=(-1)^{i+1}F_ib_{j,i}-\sum_{\ell\neq j}\sum_{r\neq i}\left((-1)^{r}\dfrac{\partial F_r}{\partial x_{j,i}}x_{\ell,r}b_{\ell,i}\right)-\sum_{r\neq i}\left((-1)^{r}\dfrac{\partial F_r}{\partial x_{j,i}}x_{j,r}b_{j,i}\right)\\
    &=(-1)^{i+1}F_ib_{j,i}-\sum_{\ell\neq j}\sum_{r\neq i}\left((-1)^{r}\dfrac{\partial F_r}{\partial x_{j,i}}x_{\ell,r}b_{\ell,i}\right)-(-1)^{i+1}F_ib_{j,i}\tag{3}\\
    &=(-1)^{i+1}F_ib_{j,i}-0-(-1)^{i+1}F_ib_{j,i}\tag{3}\\
    &=0.
    \end{align*}

    It follows that we may take $e_i \cdot b_{j,i}=0$.

\textbf{Case:} $j=1$, $i \neq k$.  

When $j=1$, we claim that
\[
e_ib_{1,k}=(-1)^{i+1}\sum_{\ell=2}^n \dfrac{\partial F_i}{\partial x_{\ell,k}}c_{\ell,1}.
\]

We may take $i \neq k$.  Consider
\begin{align*}
\partial(e_ib_{1,k})&=(-1)^{i+1}F_ib_{1,k}-\sum_{r \neq i}\left(\dfrac{\partial F_r}{\partial x_{1,k}}(-1)^r\left(\sum_{\ell=1}^nx_{\ell,r}b_{\ell,i}\right)\right)-\dfrac{\partial F_i}{\partial x_{1,k}}(-1)^{i+1}\left(\sum_{t \neq i}x_{1,t}b_{1,t}-\sum_{u \neq 1}x_{u,i}b_{u,i}\right)\\
&=(-1)^{i+1}F_ib_{1,k}-\sum_{r \neq i}\left(\dfrac{\partial F_r}{\partial x_{1,k}}(-1)^r\left(\sum_{\ell=1}^nx_{\ell,r}b_{\ell,i}\right)\right)\\
&\hspace{25.5mm}-(-1)^{i+1}\dfrac{\partial F_i}{\partial x_{1,k}}b_{1,k}+\dfrac{\partial F_i}{\partial x_{1,k}}(-1)^{i+1}\left(\sum_{t \neq i,k}x_{1,t}b_{1,t}-\sum_{u \neq 1}x_{u,i}b_{u,i}\right).  
\end{align*}

We first examine the coefficients of the basis elements $b_{m,n}$ in $\partial(e_ib_{1,k})$.  

The coefficient of $b_{\ell,i}$ for $\ell \neq 1$ is
\begin{align*}
    -\sum_{r \neq i}(-1)^r \dfrac{\partial F_r}{\partial x_{1,k}}x_{\ell,r}+(-1)^{i+1}\dfrac{\partial F_i}{\partial x_{1,k}}x_{\ell,i}=-\left(\sum_{r=1}^{n+1}(-1)^r \dfrac{\partial F_r}{\partial x_{1,k}}x_{\ell,r}\right)=0,  
\end{align*}
by Lemma \ref{syzarg}, equation (3), and Remark \ref{trivial}.  

Thus, the coefficient of $b_{\ell,i}$ for $\ell \neq 1$ in $\partial(e_ib_{1,k})$ is zero.  

The coefficient of $b_{1,i}$ is 
\[
-\sum_{r \neq i} \dfrac{\partial F_r}{\partial x_{1,k}}(-1)^rx_{1,r}.
\]

The coefficient of $b_{1,k}$ is
\[
(-1)^{i+1}F_i-(-1)^{i+1}\dfrac{\partial F_i}{\partial x_{1,k}}x_{1,k}.  
\]

The coefficient of $b_{1,t}$ for $t \neq i,k$ is
\[
\dfrac{\partial F_i}{\partial x_{1,k}}(-1)^{i}x_{1,t}. 
\]

We now examine the coefficients of the basis elements $b_{m,n}$ in $\partial\left((-1)^{i+1}\sum_{\ell=2}^n \dfrac{\partial F_i}{\partial x_{\ell,k}}c_{\ell,1}\right)$.  We have
\begin{align*}
\partial\left((-1)^{i+1}\sum_{\ell=2}^n \dfrac{\partial F_i}{\partial x_{\ell,k}}c_{\ell,1}\right)&=(-1)^{i+1}\sum_{\ell=2}^n\left(\dfrac{\partial F_i}{\partial x_{\ell,k}}\left(\sum_{s=1}^{n+1}x_{\ell,s}b_{1,s}\right)\right)\\
&=(-1)^{i+1}\sum_{\ell=2}^n\sum_{s=1}^{n+1}\dfrac{\partial F_i}{\partial x_{\ell,k}}x_{\ell,s}b_{1,s}\\
&=(-1)^{i+1}\sum_{\ell=2}^n\sum_{s \neq k}\dfrac{\partial F_i}{\partial x_{\ell,k}}x_{\ell,s}b_{1,s}+(-1)^{i+1}\sum_{\ell=2}^n\dfrac{\partial F_i}{\partial x_{\ell,k}}x_{\ell,k}b_{1,k}\\
&=(-1)^{i+1}\sum_{\ell=2}^n\sum_{s \neq k}\dfrac{\partial F_i}{\partial x_{\ell,k}}x_{\ell,s}b_{1,s}+(-1)^{i+1}\left(F_i-\dfrac{\partial F_i}{\partial x_{1,k}}x_{1,k}\right)b_{1,k}.
\end{align*}

The coefficient of $b_{\ell,i}$ for $\ell \neq 1$ is zero.  

The coefficient of $b_{1,i}$ is
\[
(-1)^{i+1}\sum_{\ell=2}^n \dfrac{\partial F_i}{\partial x_{\ell,k}}x_{\ell,i}
\]

We claim that
\[
(-1)^{i+1}\sum_{\ell=2}^n \dfrac{\partial F_i}{\partial x_{\ell,k}}x_{\ell,i}=-\sum_{r \neq i} \dfrac{\partial F_r}{\partial x_{1,k}}(-1)^rx_{1,r}.
\]
By Lemma \ref{syzarg}, equation (2),
\begin{align*}
(-1)^{i+1}\sum_{\ell=2}^n \dfrac{\partial F_i}{\partial x_{\ell,k}}x_{\ell,i}&=(-1)^{i+1}\left((-1)^{i+k+1}F_k-\dfrac{\partial F_i}{\partial x_{1,k}}x_{1,i}\right)\\
&=(-1)^kF_k+(-1)^i\dfrac{\partial F_i}{\partial x_{1,k}}x_{1,i}.
\end{align*}
Also,
\begin{align*}
-\sum_{r \neq i} \dfrac{\partial F_r}{\partial x_{1,k}}(-1)^rx_{1,r}&=-\sum_{r \neq i,k} \dfrac{\partial F_r}{\partial x_{1,k}}(-1)^rx_{1,r}\\
&=-\sum_{r \neq k} \dfrac{\partial F_r}{\partial x_{1,k}}(-1)^rx_{1,r}+(-1)^i\dfrac{\partial F_i}{\partial x_{1,k}}x_{1,i}\\
&=(-1)^{k+1+1}F_k+(-1)^i\dfrac{\partial F_i}{\partial x_{1,k}}x_{1,i}\tag{3}\\
&=(-1)^{k}F_k+(-1)^i\dfrac{\partial F_i}{\partial x_{1,k}}x_{1,i}.
\end{align*}

The coefficient of $b_{1,k}$ is
\[
(-1)^{i+1}\left(F_i-\dfrac{\partial F_i}{\partial x_{1,k}}x_{1,k}\right).
\]

The coefficient of $b_{1,t}$ for $t \neq i,k$ is
\[
(-1)^{i+1}\sum_{\ell=2}^n \dfrac{\partial F_i}{\partial x_{\ell,k}}x_{\ell,t}
\]

By Lemma \ref{lem}, equation (4), 
\[
(-1)^{i+1}\sum_{\ell=2}^n \dfrac{\partial F_i}{\partial x_{\ell,k}}x_{\ell,t}=\dfrac{\partial F_i}{\partial x_{1,k}}(-1)^{i}x_{1,t}.
\]

Since all coefficients of the basis elements $b_{m,n}$ match, we may take
\[
e_ib_{1,k}=(-1)^{i+1}\sum_{\ell=2}^n \dfrac{\partial F_i}{\partial x_{\ell,k}}c_{\ell,1}.
\]


\textbf{Case:} $j \neq 1$, $i \neq k$.

Assume $j \neq 1$, and take $i \neq k$.  We claim that
\[
e_ib_{j,k}=(-1)^{i+1}\sum_{\ell=1}^n \dfrac{\partial F_i}{\partial x_{\ell,k}}c_{\ell,j}.
\]

Consider
\begin{align*}
\partial(e_ib_{j,k})&=(-1)^{i+1}F_ib_{j,k}-\sum_{r \neq i}\left(\dfrac{\partial F_r}{\partial x_{j,k}}(-1)^r\left(\sum_{\ell=1}^nx_{\ell,r}b_{\ell,i}\right)\right)-\dfrac{\partial F_i}{\partial x_{j,k}}(-1)^{i+1}\left(\sum_{t \neq i}x_{1,t}b_{1,t}-\sum_{u \neq 1}x_{u,i}b_{u,i}\right)\\
&=(-1)^{i+1}F_ib_{j,k}-\sum_{r \neq i}\sum_{\ell=1}^n(-1)^r\dfrac{\partial F_r}{\partial x_{j,k}}x_{\ell,r}b_{\ell,i}-\dfrac{\partial F_i}{\partial x_{j,k}}(-1)^{i+1}\left(\sum_{t \neq i}x_{1,t}b_{1,t}-\sum_{u \neq 1}x_{u,i}b_{u,i}\right)\\
&=(-1)^{i+1}F_ib_{j,k}-\sum_{r \neq i}\left(\sum_{\ell=1}^{j-1}(-1)^r\dfrac{\partial F_r}{\partial x_{j,k}}x_{\ell,r}b_{\ell,i}+\sum_{\ell=j+1}^n(-1)^r\dfrac{\partial F_r}{\partial x_{j,k}}x_{\ell,r}b_{\ell,i}\right)\\
\hspace{.25mm}&-\sum_{r \neq i}(-1)^r\dfrac{\partial F_r}{\partial x_{j,k}}x_{j,r}b_{j,i}-\dfrac{\partial F_i}{\partial x_{j,k}}(-1)^{i+1}\left(\sum_{t \neq i}x_{1,t}b_{1,t}-\sum_{u \neq 1}x_{u,i}b_{u,i}\right)\\
&=(-1)^{i+1}F_ib_{j,k}-\sum_{r \neq i}\left(\sum_{\ell=1}^{j-1}(-1)^r\dfrac{\partial F_r}{\partial x_{j,k}}x_{\ell,r}b_{\ell,i}+\sum_{\ell=j+1}^n(-1)^r\dfrac{\partial F_r}{\partial x_{j,k}}x_{\ell,r}b_{\ell,i}\right)\tag{3}\\
\hspace{.25mm}&-\left((-1)^{k+1}F_k-(-1)^i\dfrac{\partial F_i}{\partial x_{j,i}}\right)b_{j,i}-\dfrac{\partial F_i}{\partial x_{j,k}}(-1)^{i+1}\left(\sum_{t \neq i}x_{1,t}b_{1,t}-\sum_{u \neq 1}x_{u,i}b_{u,i}\right).
\end{align*}

We first examine the coefficients of the basis elements $b_{m,n}$ in $\partial(e_ib_{j,k})$.  

The coefficient of $b_{j,k}$ is
\[
(-1)^{i+1}F_i.
\]

The coefficient of $b_{1,v}$ for $v \neq i$ is
\[
(-1)^i\dfrac{\partial F_i}{\partial x_{j,k}}x_{1,v}.
\]

The coefficient in $b_{1,i}$ is
\[
-\sum_{r \neq i}\dfrac{\partial F_r}{\partial x_{j,k}}(-1)^rx_{1,r}.
\]

When $v=i$, we claim no discrepancy, i.e.,
\[
(-1)^i\dfrac{\partial F_i}{\partial x_{j,k}}x_{1,i}=-\sum_{r \neq i}\dfrac{\partial F_r}{\partial x_{j,k}}(-1)^rx_{1,r}.
\]
To see this, let $Y$ be the matrix obtained from $X_k$ by replacing the row indexed by $j$ by a second row indexed by $1$.  The determinant of $Y$ is zero, and one way to compute $\det Y$ is via Laplace expansion along the first row:
\begin{align*}
0=\det Y&=\pm \left(\sum_{r < k}(-1)^{1+r}x_{1,r}\det X_{k,r}^1+\sum_{r > k}(-1)^{1+r+1}x_{1,r}\det X_{k,r}^1\right)\\
&=\pm\left(\sum_{r <k}(-1)^{1+r}x_{1,r}(-1)^{1+k+1}\dfrac{\partial F_r}{\partial x_{1,k}}+\sum_{r>k}(-1)^{1+r+1}x_{1,r}(-1)^{1+k}\dfrac{\partial F_r}{\partial x_{1,k}}\right)\tag{Remark \ref{swap}}\\
&=\pm(-1)^{k+1}\sum_{r=1}^{n+1}(-1)^rx_{1,r}\dfrac{\partial F_r}{\partial x_{j,k}}.
\end{align*}

The coefficient of $b_{u,i}$ for $u \neq 1$ is
\begin{align*}
-\sum_{r \neq i}\left(\dfrac{\partial F_r}{\partial x_{j,k}}(-1)^rx_{u,r}\right)+(-1)^{i+1}\dfrac{\partial F_i}{\partial x_{j,k}}x_{u,i}=-\sum_{r=1}^{n+1}\dfrac{\partial F_r}{\partial x_{j,k}}(-1)^rx_{u,r}.
\end{align*}

When $j=u$, Lemma \ref{syzarg}, equation (3), tell us that 
\[
-\sum_{r=1}^{n+1}\dfrac{\partial F_r}{\partial x_{j,k}}(-1)^rx_{u,r}=(-1)^kF_k.
\]

When $j \neq u$,  
\[
\sum_{r=1}^{n+1}\dfrac{\partial F_r}{\partial x_{j,k}}(-1)^rx_{u,r}=0,
\]
by Lemma \ref{syzarg}, equation (3), and Remark \ref{trivial}.  

Thus, the coefficient of $b_{u,i}$ for $u \neq 1$ and $j \neq u$ is zero.  

The coefficients of all other $b_{m,n}$ are zero.  

We now examine the coefficients of the basis elements $b_{m,n}$ in $\partial\left((-1)^{i+1}\sum_{\ell=1}^n \dfrac{\partial F_i}{\partial x_{\ell,k}}c_{\ell,j}\right)$.  We have
\begin{align*}
\partial\left((-1)^{i+1}\sum_{\ell=1}^n \dfrac{\partial F_i}{\partial x_{\ell,k}}c_{\ell,j}\right)&=(-1)^{i+1}\sum_{\ell=1}^n\left(\dfrac{\partial F_i}{\partial x_{\ell,k}}\left(\sum_{u=1}^{n+1}x_{\ell,u}b_{j,u}-\delta_{\ell,j}\sum_{v=1}^{n+1}x_{1,v}b_{1,v}\right)\right)\\
&=(-1)^{i+1}\left(\sum_{\ell=1}^n\sum_{u=1}^{n+1}\dfrac{\partial F_i}{\partial x_{\ell,k}}x_{\ell,u}b_{j,u}-\dfrac{\partial F_i}{\partial x_{j,k}}\sum_{v=1}^{n+1}x_{1,v}b_{1,v}\right)\\
&=(-1)^{i+1}\left(\sum_{\ell=1}^n\dfrac{\partial F_i}{\partial x_{\ell,k}}x_{\ell,k}b_{j,k}+\sum_{\ell=1}^n\dfrac{\partial F_i}{\partial x_{\ell,k}}x_{\ell,i}b_{j,i}-\dfrac{\partial F_i}{\partial x_{j,k}}\sum_{v=1}^{n+1}x_{1,v}b_{1,v}\right)\tag{4}\\ 
&=(-1)^{i+1}\left(F_ib_{j,k}+\sum_{\ell=1}^n\dfrac{\partial F_i}{\partial x_{\ell,k}}x_{\ell,i}b_{j,i}-\dfrac{\partial F_i}{\partial x_{j,k}}\sum_{v=1}^{n+1}x_{1,v}b_{1,v}\right)\tag{1}.
\end{align*}

The coefficient of $b_{j,k}$ is $(-1)^{i+1}F_i$.

The coefficient of $b_{1,v}$ is
\[
(-1)^i\dfrac{\partial F_i}{\partial x_{j,k}}x_{1,v}.
\]

The coefficient of $b_{j,i}$ is
\[
(-1)^{i+1}\sum_{\ell=1}^n\dfrac{\partial F_i}{\partial x_{\ell,k}}x_{\ell,i},
\]
which, by Lemma \ref{syzarg}, equation 2, is equal to
\[
(-1)^{i+1}(1)^{i+k+1}F_k=(-1)^kF_k.
\]

The coefficients of all other $b_{m,n}$ are zero.  

\textbf{End $e_ib_{j,k}$ action.}

\textbf{Begin $T_ia_j$ action:}
    
    We claim that 
    \[
    T_i\cdot a_j=(-1)^j\sum_{k=1}^nx_{k,j}c_{i,k}.
    \]

    We have
    \begin{align*}
    \partial(T_ia_j)&=\sum_{\ell=1}^{n+1}x_{i,\ell}e_{\ell}a_j\\
    &=\sum_{\ell \neq j}\left(x_{i,\ell}\left((-1)^j\sum_{k=1}^n x_{k,j}b_{k,\ell}\right)\right)+(-1)^{j+1}\left(x_{i,j}\sum_{u \neq j}x_{1,u}b_{1,u}-\sum_{v \neq 1}x_{v,j}b_{v,j}\right)\\
    &=\sum_{\ell \neq j}\left(x_{i,\ell}\left((-1)^j\sum_{k=1}^n x_{k,j}b_{k,\ell}\right)\right)+(-1)^{j+1}x_{i,j}\sum_{u=1}^{n+1}x_{1,u}b_{1,u}\\&+(-1)^{j+1}x_{i,j}x_{1,j}b_{1,j}-(-1)^{j+1}x_{i,j}x_{1,j}b_{1,j}-(-1)^{j+1}x_{i,j}\sum_{v=1}^nx_{v,j}b_{v,j}\\
    &=\sum_{\ell \neq j}\left(x_{i,\ell}\left((-1)^j\sum_{k=1}^n x_{k,j}b_{k,\ell}\right)\right)+(-1)^{j+1}x_{i,j}\sum_{u=1}^{n+1}x_{1,u}b_{1,u}-(-1)^{j+1}x_{i,j}\sum_{v=1}^nx_{v,j}b_{v,j}.
    \end{align*}

    Also,
    \begin{align*}
        \partial\left((-1)^j\sum_{k=1}^nx_{k,j}c_{i,k}\right)&=(-1)^j\sum_{k=1}^n\left(x_{k,j}\left(\sum_{\ell=1}^{n+1}x_{i,\ell}b_{k,\ell}-\delta_{i,k}\sum_{v=1}^{n+1}x_{1,v}b_{1,v}\right)\right)\\
        &=(-1)^j\sum_{k\neq i}\sum_{\ell=1}^{n+1}x_{k,j}x_{i,\ell}b_{k,\ell}+(-1)^jx_{i,j}\left(\sum_{\ell=1}^{n+1}x_{i,\ell}b_{i,\ell}-\sum_{v=1}^{n+1}x_{1,v}b_{1,v}\right)\\
        &=(-1)^j\sum_{\ell=1}^{n+1}\left(x_{i,\ell}\left(\sum_{k=1}^nx_{k,j}b_{k,\ell}-x_{i,j}b_{i,\ell}\right)\right)+(-1)^jx_{i,j}\left(\sum_{\ell=1}^{n+1}x_{i,\ell}b_{i,\ell}-\sum_{v=1}^{n+1}x_{1,v}b_{1,v}\right)\\
        &=(-1)^j\sum_{\ell=1}^{n+1}\left(x_{i,\ell}\sum_{k=1}^nx_{k,j}b_{k,\ell}\right)-(-1)^j\sum_{\ell=1}^{n+1}x_{i,\ell}x_{i,j}b_{i,\ell}\\&+(-1)^j\sum_{\ell=1}^{n+1}x_{i,j}x_{i,\ell}b_{i,\ell}-(-1)^j\sum_{v=1}^{n+1}x_{i,j}x_{1,v}b_{1,v}\\
        &=(-1)^j\sum_{\ell\neq j}\left(x_{i,\ell}\sum_{k=1}^nx_{k,j}b_{k,\ell}\right)+(-1)^jx_{i,j}\sum_{k=1}^nx_{k,j}b_{k,j}-(-1)^j\sum_{\ell=1}^{n+1}x_{i,\ell}x_{i,j}b_{i,\ell}\\&+(-1)^j\sum_{\ell=1}^{n+1}x_{i,j}x_{i,\ell}b_{i,\ell}-(-1)^j\sum_{v=1}^{n+1}x_{i,j}x_{1,v}b_{1,v}\\
        &=(-1)^j\sum_{\ell\neq j}\left(x_{i,\ell}\sum_{k=1}^nx_{k,j}b_{k,\ell}\right)+(-1)^{j+1}x_{i,j}\sum_{v=1}^{n+1}x_{1,v}b_{1,v}-(-1)^{j+1}x_{i,j}\sum_{k=1}^nx_{k,j}b_{k,j},
        \end{align*}
    which is $\partial(T_ia_j)$.

    Thus,
    \[
    T_i\cdot a_j=(-1)^j\sum_{k=1}^nx_{k,j}c_{i,k}.
    \]
    
\textbf{End $T_ia_j$ action.}

It remains to show that the multiplication in the dg structure we have defined is associative.  For degree reasons, most triple products are zero.  It is not yet clear that $(e_ie_j)a_k=e_i(e_ja_k)$. 

\textbf{Case:} $i=j$.

When $i=j$, $(e_ie_j)a_k=(e_ie_i)a_k=0$.  If we further assume $i \neq k$, then
\begin{align*}
e_i(e_ia_k)&=e_i\left((-1)^k\sum_{u=1}^nx_{u,k}b_{u,i}\right)\\
&=(-1)^k\sum_{u=1}^nx_{u,k}e_ib_{u,i}\\
&=0.
\end{align*}

If instead we assume $i=k$, then
\begin{align*}
e_i(e_ia_k)&=e_i\left((-1)^{i+1}\left(\sum_{u \neq i}x_{1,u}b_{1,u}-\sum_{v \neq 1}x_{v,i}b_{v,i}\right)\right)\\
&=(-1)^{i+1}\sum_{u \neq i}x_{1,u}e_ib_{1,u}\\
&=(-1)^{i+1}\sum_{u=1}^{n+1}\left(x_{1,u}\left((-1)^{i+1}\sum_{\ell=1}^n\dfrac{\partial F_i}{\partial x_{\ell,u}}c_{\ell,1}\right)\right)\\
&=\sum_{u=1}^{n+1}\sum_{\ell=1}^nx_{1,u}\dfrac{\partial F_i}{\partial x_{\ell,u}}c_{\ell,1}.
\end{align*}
The coefficient of $c_{\ell,1}$ for $\ell \neq 1$ above is
\[
\sum_{u=1}^{n+1}x_{1,u}\dfrac{\partial F_i}{\partial x_{\ell,u}}=0,
\]
by Lemma \ref{syzarg}, equation (1).  Thus,
\[
e_i(e_ia_k)=0.  
\]


\textbf{Case:} $i \neq j$.  

When $i<j$,
\begin{align*}
(e_ie_j)a_k&=\sum_{u=1}^n(-1)^{i+j+u}\det X^u_{i,j}T_ua_k\\
&=\sum_{u=1}^n\left((-1)^{i+j+u}\det X^u_{i,j}(-1)^k\sum_{v=1}^nx_{v,k}c_{u,v}\right)\\
&=(-1)^{k}\sum_{v=1}^n\sum_{u=1}^n(-1)^{i+j+u}\det X^u_{i,j}x_{v,k}c_{u,v}\\
&=(-1)^{k}\sum_{v=1}^n\sum_{u=1}^n(-1)^{i+1}\dfrac{\partial F_i}{\partial x_{u,j}}x_{v,k}c_{u,v}. \tag{a}
\end{align*}

When $i>j$,
\begin{align*}
(e_ie_j)a_k&=-\sum_{u=1}^n(-1)^{i+j+u}\det X^u_{i,j}T_ua_k\\
&=-\sum_{u=1}^n\left((-1)^{i+j+u}\det X^u_{i,j}(-1)^{k+1}\sum_{v=1}^nx_{v,k}c_{u,v}\right)\\
&=(-1)^{k+1}\sum_{v=1}^n\sum_{u=1}^n(-1)^{i+j+u}\det X^u_{i,j}x_{v,k}c_{u,v}\\
&=(-1)^{k+1}\sum_{v=1}^n\sum_{u=1}^n(-1)^{i}\dfrac{\partial F_i}{\partial x_{u,j}}x_{v,k}c_{u,v}. \tag{b}
\end{align*}

If we further assume that $j \neq k$, then
\begin{align*}
e_i(e_ja_k)&=e_i\left((-1)^{k}\displaystyle\sum_{v=1}^nx_{v,k}b_{v,j}\right)\\
&=(-1)^{k}\displaystyle\sum_{v=1}^nx_{v,k}e_ib_{v,j}\\
&=(-1)^{k}\displaystyle\sum_{v=1}^n\left(x_{v,k}\left((-1)^{i+1}\displaystyle\sum_{u=1}^n\dfrac{\partial F_i}{\partial x_{u,j}}c_{u,v}\right)\right)\\
&=(-1)^{i+k+1}\displaystyle\sum_{v=1}^n\sum_{u=1}^nx_{v,k}\dfrac{\partial F_i}{\partial x_{u,j}}c_{u,v}.  \tag{c}
\end{align*}

If instead we assume $j=k$, then
\begin{align*}
e_i(e_ja_k)=e_i(e_ja_j)&=e_i\left((-1)^{j+1}\left(\sum_{u \neq j}x_{1,u}b_{1,u}-\sum_{v \neq 1}x_{v,j}b_{v,j}\right)\right)\\
&=(-1)^{j+1}\left(\sum_{u \neq j}x_{1,u}e_ib_{1,u}-\sum_{v \neq 1}x_{v,j}e_ib_{v,j}\right)\\
&=(-1)^{j+1}\left(\sum_{u \neq j}x_{1,u}\left((-1)^{i+1}\sum_{\ell=1}^n\dfrac{\partial F_i}{\partial x_{\ell,u}}c_{\ell,1}\right)-\sum_{v \neq 1}x_{v,j}\left((-1)^{i+1}\sum_{s=1}^n\dfrac{\partial F_i}{\partial x_{s,j}}c_{s,v}\right)\right)\\
&=(-1)^{j+1}\left((-1)^{i+1}\sum_{u \neq j}\sum_{\ell=1}^nx_{1,u}\dfrac{\partial F_i}{\partial x_{\ell,u}}c_{\ell,1}-(-1)^{i+1}\sum_{v \neq 1}\sum_{s=1}^nx_{v,j}\dfrac{\partial F_i}{\partial x_{s,j}}c_{s,v}\right)\\
&=(-1)^{j+1}\left((-1)^{i}\sum_{\ell=1}^nx_{1,j}\dfrac{\partial F_i}{\partial x_{\ell,j}}c_{\ell,1}-(-1)^{i+1}\sum_{v \neq 1}\sum_{s=1}^nx_{v,j}\dfrac{\partial F_i}{\partial x_{s,j}}c_{s,v}\right) \\
&=(-1)^{i+j+1}\left(\sum_{\ell=1}^nx_{1,j}\dfrac{\partial F_i}{\partial x_{\ell,j}}c_{\ell,1}+\sum_{v \neq 1}\sum_{s=1}^nx_{v,j}\dfrac{\partial F_i}{\partial x_{s,j}}c_{s,v}\right)\\
&=(-1)^{i+j+1}\sum_{v=1}^n\sum_{s=1}^nx_{v,j}\dfrac{\partial F_i}{\partial x_{s,j}}c_{s,v}\\
&=(-1)^{i+k+1}\sum_{v=1}^n\sum_{s=1}^nx_{v,k}\dfrac{\partial F_i}{\partial x_{s,j}}c_{s,v}. \tag{d}
\end{align*}

It remains to note that (a), (b), (c), and (d) are all equal.  


\end{proof}

\begin{cor}\label{golodsyz}
The $R$-module $\Der_{R \mid k}$ is Golod.  
\end{cor}

\begin{proof}
     By Lemma \ref{act}, the action of $\mathcal{A}_{\geq 1}$ on $\mathcal{U}$ lands in $\m_Q\mathcal{U}$, and so applying Proposition \ref{golodcon} yields that the $R$-module $\coker J^T_R$ is Golod.  By a theorem of Levin (\cite{golodsyz}, Theorem 1.1), the first syzygy of a Golod module is Golod.    
\end{proof}

\section{The Resolution of $\Der_{R \mid k}$}\label{resder}

Having established the ingredients needed to produce the relative bar resolution of $\coker J^T_R$ over $R$, we do so after introducing notation for the derivations coming from the dg action defined in the previous section.  Truncating, we obtain the minimal $R$-free resolution of $\Der_{R \mid k}$.    

\begin{notation}
Let $L_{i,j}$ be the derivation corresponding to $e_i \cdot a_j$, i.e.,
    \[
    L_{i,j}=\begin{cases}
    (-1)^j\displaystyle\sum_{k=1}^nx_{k,j}\dfrac{\partial}{\partial x_{k,i}},    & i \neq j,\\\\
     (-1)^{i+1}\left(\displaystyle\sum_{k \neq i}x_{1,k}\dfrac{\partial}{\partial x_{1,k}}-\displaystyle\sum_{\ell \neq 1}x_{\ell,i}\dfrac{\partial}{\partial x_{\ell,i}}\right),   & i=j.
    \end{cases}
    \]
    Let $M$ be the $n(n+1) \times (n+1)^2$ matrix with rows indexed by
    \[
    b_{1,1},\,b_{1,2},\ldots,b_{1,n+1},\,b_{2,1},\ldots,b_{n+1,n+1}
    \]
    and columns corresponding (via $\frac{\partial}{\partial x_{i,j}} \leftrightarrow b_{i,j}$) to
    \[
    L_{1,1},\,L_{1,2},\ldots,L_{1,n+1},\,L_{2,1},\ldots,L_{n+1,n+1}.
    \]

Set $M=M_{1,0}$, and define $M_{i,j}$ similarly, as the matrix coming from the multiplication $\mathcal{A}_i \cdot \mathcal{U}_j$.
\end{notation}

We now combine the ingredients from the previous sections to produce the minimal $R$-free resolution of $\coker J^T_R$.  This is done using the construction of the relative bar resolution via dg structures due to Iyengar \cite{Iyengar}.  For us, the construction requires a dg algebra resolution of $R$ over $Q$ (the Hilbert-Burch complex $\mathcal{A}$ from Proposition \ref{EN}) and a dg $\mathcal{A}$-module resolution of $\coker J^T_R$ over $Q$ (the resolution $\mathcal{U}$ from Theorem \ref{mfrcokerjt}).  We saw that $R$ is a Golod ring (Corollary \ref{Rgol}) and that $\coker J^T_R$ is a Golod module (see the proof of Corollary \ref{golodsyz}).  The Golod property together with the fact that $\mathcal{A}$ and $\mathcal{U}$ are minimal implies that the resulting relative bar resolution of $\coker J^T_R$ over $R$ is minimal.     
    
\begin{prop}\label{cokresn}
    The minimal $R$-free resolution of $\coker J^T_R$ is the relative bar resolution
    \[
    \dots \longrightarrow R^{n(n+1)(n+2)} \xrightarrow{\small\begin{bmatrix}
        J^T_R & \varphi_2\\
        M_{1,1} & M_{2,0}
    \end{bmatrix}} R^{2n(n+1)}\xrightarrow{\small\begin{bmatrix}
        \partial_2 & M_{1,0}
    \end{bmatrix}}\normalsize R^{n(n+1)}\xrightarrow{J^T_R} R^{n+1}\longrightarrow \coker J^T_R \longrightarrow 0.
    \]

\end{prop}

We recall that $\Der_{R \mid k}$ is isomorphic to the kernel of $J^T_R$.  It follows that by truncating the relative bar resolution from Proposition \ref{cokresn}, we get the minimal $R$-free resolution of $\Der_{R \mid k}$.   

\begin{thm}\label{mainresult}
    The minimal $R$-free resolution of $\Der_{R \mid k}$ is the truncated relative bar resolution
    \[
    \dots \longrightarrow R^{n(n+1)(n+2)} \xrightarrow{\small\begin{bmatrix}
        J^T_R & \varphi_2\\
        M_{1,1} & M_{2,0}
    \end{bmatrix}} R^{2n(n+1)}\xrightarrow{\small\begin{bmatrix}
        \partial_2 & M_{1,0}
    \end{bmatrix}} \Der_{R \mid k} \longrightarrow 0.
    \]
\end{thm}

Theorem \ref{mainresult} recovers a minimal generating set for $\Der_{R \mid k}$ (see \cite{bv} and \cite{CB} for other minimal generating sets of $\Der_{R \mid k}$, obtained using different approaches).  

\begin{cor}\label{mg}
    A minimal generating set for $\Der_{R \mid k}$ over $R$ is 
    \[
    \mathcal{M}:=\left\{V_{r,s},B_\ell,L_{i,j}\right\},
    \]
    where $r,s, \in \{1,\ldots,n\}$, $r \neq s$, $\ell \in \{2,\ldots,n\}$, and $i,j \in \{1,\ldots,n+1\}$. 
\end{cor}

For the sake of visualization, we give an example of what the first differential in the resolution from Theorem \ref{mainresult} looks like in unabbreviated form.  

\begin{ex}\textup{
    When $n=2$, the matrix $\begin{bmatrix}
        \partial_2 & M_{1,0}
    \end{bmatrix}$ is equal to
    \[
    \left[\begin{tabular}{cccccccccccc}
    $x_{2,1}$&$0$&$-x_{1,1}$&$0$&$x_{1,2}$&$-x_{1,3}$&$0$&$-x_{1,1}$&$0$&$0$&$0$&$x_{1,1}$\\
    $x_{2,2}$&$0$&$-x_{1,2}$&$x_{1,2}$&$0$&$0$&$-x_{1,1}$&$0$&$-x_{1,3}$&$0$&$0$&$x_{1,2}$\\
    $x_{2,3}$&$0$&$-x_{1,3}$&$x_{1,3}$&$0$&$0$&$0$&$-x_{1,3}$&$0$&$-x_{1,1}$&$x_{1,2}$&$0$\\
    $0$&$x_{1,1}$&$x_{2,1}$&$-x_{2,1}$&$x_{2,2}$&$-x_{2,3}$&$0$&$0$&$0$&$0$&$0$&$0$\\
    $0$&$x_{1,2}$&$x_{2,2}$&$0$&$0$&$0$&$-x_{2,1}$&$x_{2,2}$&$-x_{2,3}$&$0$&$0$&$0$\\
    $0$&$x_{1,3}$&$x_{2,3}$&$0$&$0$&$0$&$0$&$0$&$0$&$-x_{2,1}$&$x_{2,2}$&$-x_{2,3}$\\
    \end{tabular}\right],
    \]
    in which the columns are indexed
    \[  c_{2,1},\,c_{1,2},\,c_{2,2},\,L_{1,1},\,L_{1,2},\,L_{1,3},\,L_{2,1},\,L_{2,2},\,L_{2,3},\,L_{3,1},\,L_{3,2},\,L_{3,3}.  
    \]
    }
\end{ex} 

\begin{rmk}
    One may wonder how the Euler derivation (Example \ref{eu}) arises from the minimal set of generators from Corollary \ref{mg}.  In terms of the elements of $\mathcal{M}$, the Euler derivation $E$ can be expressed as
    \[
    E=n\left(\sum_{k=1}^{n+1}(-1)^{k+1}L_{k,k}\right)+(n+1)\left(\sum_{\ell=2}^n B_\ell\right).
    \]
    To see this, examine the coefficients of $\dfrac{\partial}{\partial x_{u,v}}$ in
    \begin{align*}
        n\left(\sum_{k=1}^{n+1}(-1)^{k+1}L_{k,k}\right)+(n+1)\left(\sum_{\ell=2}^n B_\ell\right)&=n\left(\sum_{k=1}^{n+1}(-1)^{k+1}(-1)^{k+1}\left(\sum_{u \neq k}x_{1,u}\dfrac{\partial}{\partial x_{1,u}}-\sum_{v \neq 1}x_{v,k}\dfrac{\partial}{\partial x_{v,k}}\right)\right)\\
        &+(n+1)\left(\sum_{\ell=2}^n\left(\sum_{s=1}^{n+1}x_{\ell,s}\dfrac{\partial}{\partial x_{\ell,s}}-\sum_{t=1}^{n+1}x_{1,t}\dfrac{\partial}{\partial x_{1,t}}\right)\right)\\
        &=n\left(\sum_{k=1}^{n+1}\left(\sum_{u \neq k}x_{1,u}\dfrac{\partial}{\partial x_{1,u}}-\sum_{v \neq 1}x_{v,k}\dfrac{\partial}{\partial x_{v,k}}\right)\right)\\
        &+(n+1)\left(\sum_{\ell=2}^n\left(\sum_{s=1}^{n+1}x_{\ell,s}\dfrac{\partial}{\partial x_{\ell,s}}-\sum_{t=1}^{n+1}x_{1,t}\dfrac{\partial}{\partial x_{1,t}}\right)\right).
    \end{align*}
    The coefficient of $\dfrac{\partial}{\partial x_{1,v}}$ is $n^2x_{1,v}-(n+1)(n-1)x_{1,v}=x_{1,v}$, and the coefficient of $\dfrac{\partial}{\partial x_{u,v}}$ for $u \neq 1$ is $-nx_{u,v}+(n+1)x_{u,v}=x_{u,v}$.  
\end{rmk}

The construction of the relative bar resolution via dg structures implies that the Poincar\'e series of $\Der_{R \mid k}$ is rational.  

\begin{cor}\label{poin}
    The Poincar\'e series of $\Der_{R \mid k}$ over $R$ is
    \[
    P^R_{\Der_{R \mid k}}(t)=\dfrac{n(n+1)(2+nt)}{1-t-nt^2}.
    \]
\end{cor}

\begin{proof}
    Since $\coker J^T_R$ is Golod,
    \begin{align*}
    P^R_{\coker J^T_R}(t)&=\dfrac{P^Q_{\coker J^T_R}(t)}{1-t\left(P^Q_R(t)-1\right)}\\
    &=\dfrac{n+1+n(n+1)t+(n-1)(n+1)t^2}{1-t(1+(n+1)t+nt^2-1)}\\
    &=\dfrac{(n+1)(1+nt+(n-1)t^2)}{1-(n+1)t^2-nt^3}\\
    &=\dfrac{(n+1)(1+t)(1+(n-1)t)}{(1+t)(1-t-nt^2)}\\
    &=\dfrac{(n+1)(1+(n-1)t)}{1-t-nt^2}.
    \end{align*}
    By the exact sequence
    \[
0 \to \Der_{R \mid k} \to R^{n(n+1)} \xrightarrow{J^T_R} R^{n+1} \to \coker J^T_R \to 0, 
    \]
    we have that $P^R_{\Der_{R \mid k}}(t)$ is obtained by subtracting $(n+1)+n(n+1)t$ and multiplying by $\dfrac{1}{t^2}$:
    \begin{align*}
    P^R_{\Der_{R \mid k}}(t)=\left(\dfrac{2nt^2+2n^2t^2+n^2t^3+n^3t^3}{1-t-nt^2}\right)\cdot\dfrac{1}{t^2}=\dfrac{n(n+1)(2+nt)}{1-t-nt^2}.
    \end{align*}
    \end{proof}

As a final corollary, we note the linearity of the resolution of $\Der_{R \mid k}$ when $n=2$.  

\begin{cor}
    When $n=2$, the minimal $R$-free resolution of $\Der_{R \mid k}$ is linear, i.e., the entries of the matrices representing the differentials in the minimal $R$-free resolution of $\Der_{R \mid k}$ are either zero or variables.  
\end{cor}

\begin{proof}
    This is because partials of minors in the $2 \times 3$ case are just variables.  
\end{proof}

\section*{Acknowledgments}

Special thanks to Claudia Miller for all of her guidance on this project.  Thanks also to Josh Pollitz and Jack Jeffries for helpful discussions and to the reviewers for their comments and suggestions.  Several results in this paper I first conjectured based on examples I computed using Macaulay2 \cite{M2}. 


\printbibliography

@book{bv,
title = {Determinantal Rings},
    author = {Winfried Bruns and Udo Vetter},
    isbn = {3540194681},
    series = {Lecture Notes in Mathematics},
    year = {1988},
    publisher = {Spring-Verlag}
}

@book{bh,
place={Cambridge}, edition={2}, series={Cambridge Studies in Advanced Mathematics}, title={Cohen-Macaulay Rings}, publisher={Cambridge University Press}, author={Bruns, Winfried and Herzog, Jürgen}, year={1998}, collection={Cambridge Studies in Advanced Mathematics}, isbn = {9780511608681}
}

@article{hema,
    author = "Hema Srinivasan",
    title = "Algebra Structures on Some Canonical Resolutions",
    journal = "Journal of Algebra",
    volume = "122",
    pages = "150--187",
    year = "1989",
    DOI = "http://dx.doi.org/10.1002/andp.19053221004"
}

@article{BEA,
    author = "David A. Buchsbaum and David Eisenbud",
    title = "What Makes a Complex Exact?",
    journal = "Journal of Algebra",
    volume = "25",
    pages = "259--268",
    year = "1973",
    DOI = "https://doi.org/10.1016/0021-8693(73)90044-6"
}

@InProceedings{dga,
author="Beck, Kristen A.
and Sather-Wagstaff, Sean",
editor="Cooper, Susan M.
and Sather-Wagstaff, Sean",
title="A Somewhat Gentle Introduction to Differential Graded Commutative Algebra",
booktitle="Connections Between Algebra, Combinatorics, and Geometry",
year="2014",
publisher="Springer New York",
address="New York, NY",
pages="3--99",
isbn="978-1-4939-0626-0"
}

@Inbook{IFR,
author="Avramov, Luchezar L.",
editor="Elias, J.
and Giral, J. M.
and Mir{\'o}-Roig, R. M.
and Zarzuela, S.",
title="Infinite Free Resolutions",
bookTitle="Six Lectures on Commutative Algebra",
year="1998",
publisher="Birkh{\"a}user Basel",
address="Basel",
pages="1--118",
isbn="978-3-0346-0329-4",
doi="10.1007/978-3-0346-0329-4_1",
url="https://doi.org/10.1007/978-3-0346-0329-4_1"
}

@article{Gor,
 ISSN = {03713539},
 URL = {http://www.jstor.org/stable/43698822},
 author = {Shiro Goto},
 journal = {Science Reports of the Tokyo Kyoiku Daigaku, Section A},
 number = {329/346},
 pages = {129--145},
 publisher = {Editorial Committee of Tsukuba Journal of Mathematics},
 title = {When do the determinantal ideals define Gorenstein rings?},
 volume = {12},
 year = {1974}
}

@article{Iyengar,
    author = "Srikanth Iyengar",
    title = "Free Resolutions and Change of Rings",
    journal = "Journal of Algebra",
    volume = "190",
    pages = "195--213",
    year = "1997",
    DOI = "https://doi.org/10.1006/jabr.1996.6901"
}

@misc{Burke,
      title={Higher homotopies and Golod rings}, 
      author={Jesse Burke},
      year={2015},
      eprint={1508.03782},
      archivePrefix={arXiv},
      primaryClass={math.AC},
      url={https://arxiv.org/abs/1508.03782}, 
}

@article{hilbert,
  title={Uber die Theorie der algebraischen Formen},
  author={David R. Hilbert},
  journal={Mathematische Annalen},
  year={1890},
  volume={36},
  pages={473-534},
  url={https://api.semanticscholar.org/CorpusID:179177713}
}

@article{burch, title={On ideals of finite homological dimension in local rings}, volume={64}, DOI={10.1017/S0305004100043620}, number={4}, journal={Mathematical Proceedings of the Cambridge Philosophical Society}, author={Burch, Lindsay}, year={1968}, pages={941–948}}

@article{herzog, title={Komplexe, auflosungen, und dualit\"at in der lokalen algebra}, journal={Habilitationsschrift}, author={Herzog, J\"urgen}, year={1973}}

@article{rc,
      title={Resolutions of differential operators of low order for an isolated hypersurface singularity}, 
      author={Rachel N. Diethorn and Jack Jeffries and Claudia Miller and Nicholas Packauskas and Josh Pollitz and Hamidreza Rahmati and Sophia Vassiliadou},
      year={to appear},
      journal={Michigan Mathematical Journal}
}

@article{yu,
title = {A free resolution of the module of derivations for generic arrangements},
journal = {Journal of Algebra},
volume = {136},
number = {2},
pages = {432-438},
year = {1991},
issn = {0021-8693},
doi = {https://doi.org/10.1016/0021-8693(91)90054-C},
url = {https://www.sciencedirect.com/science/article/pii/002186939190054C},
author = {Sergey Yuzvinsky}
}

@Misc{M2,
          author = {Grayson, Daniel R. and Stillman, Michael E.},
          title = {Macaulay2, a software system for research in algebraic geometry},
          howpublished = {Available at \url{http://www2.macaulay2.com}}
        }

@article{golodsyz,
title = {Modules and Golod homomorphisms},
journal = {Journal of Pure and Applied Algebra},
volume = {38},
number = {2},
pages = {299-304},
year = {1985},
issn = {0022-4049},
doi = {https://doi.org/10.1016/0022-4049(85)90017-9},
url = {https://www.sciencedirect.com/science/article/pii/0022404985900179},
author = {Gerson Levin}
}

@book{cawavtag,
  title={Commutative Algebra: with a View Toward Algebraic Geometry},
  author={David Eisenbud},
  isbn={9783540942696},
  lccn={lc94017351},
  series={Graduate Texts in Mathematics},
  year={1995},
  publisher={Springer-Verlag}
}

@article{CB,
  title={The module of logarithmic derivations of a generic determinantal ideal},
  journal={Proceedings of the American Mathematical Society},
  author={Ricardo Burity and Cleto B. Miranda-Neto},
  year={2020},
  url={https://api.semanticscholar.org/CorpusID:219092477}
}
  
\end{document}